\documentclass[12pt]{article}
\usepackage{amssymb, amsmath, amscd, amsthm, setspace, latexsym}
\begin{document}
\title{On pairs of matrices generating matrix rings and their presentations}
\author{B.V.~Petrenko\\
\small{Department of Mathematics,}\\
\small{Texas A\&M University, College Station, TX 77843-3368, USA.}\\
\small{E-mail: petrenko@math.tamu.edu}
\and
S.N.~Sidki\\
\small{Department of Mathematics,}\\
\small{University of Bras\'{\i}lia, 70.910 Bras\'{\i}lia DF, BRAZIL}\\
\small{E-mail:sidki@mat.unb.br}}

\newtheorem{theorem}{Theorem}[section]
\newtheorem{lemma}[theorem]{Lemma}
\newtheorem{proposition}[theorem]{Proposition}
\newtheorem{corollary}[theorem]{Corollary} 
\newtheorem{remark}[theorem]{Remark}
\newtheorem{conjecture}[theorem]{Conjecture}
\newtheorem{definition}[theorem]{Definition}
\newtheorem{example}[theorem]{Example}
\newtheorem{property}[theorem]{Property}
\newtheorem{question}[theorem]{Question}

\maketitle

\begin{abstract} Let $M_n(\mathbb{Z})$
the ring of $n$-by-$n$ matrices with integral entries, and 
$n \geq 2$. 
This paper studies the set $G_n(\mathbb{Z})$ of pairs 
$(A,B) \in M_n(\mathbb{Z})^2$ generating $M_n(\mathbb{Z})$ as a ring.
We use several presentations of 
$M_{n}(\mathbb{Z})$ with generators $X=\sum_{i=1}^n E_{i+1,i}$ and $Y=E_{11}$
to obtain the following consequences.
\begin{enumerate}
\item Let $k \geq 1$. Then the rings $M_n(\mathbb{Q})^k$ and  
$\bigoplus_{j=1}^{k} M_{n_j} (\mathbb{Z})$, where
$n_1, \ldots, n_k \geq 2$ are pairwise relatively prime,
have presentations with $2$ generators and finitely many relations.

\item Let $D$ be a commutative domain of sufficiently large characteristic
over which every finitely generated projective
module is free. We use
$4$ relations for $X$ and $Y$ to describe all representations of the ring 
$M_{n}(D)$ into $M_{m}(D)$ for $m \geq n$.

\item We obtain information about the asymptotic density of $G_n(F)$ in $M_n(F)^2$
over different fields, and over the integers.

\bigskip

\end{enumerate}
\end{abstract}

\tableofcontents

\section{Introduction}
\subsection{Terminology and notation}
All rings in this paper, often denoted by $R$, are assumed associative 
with a two-sided identity element, unless stated otherwise. We
denote by $\mathcal{U}(R)$ the unit group of $R$. 
We do not assume that a subring of a ring necessarily
contains the identity element of the ring.
All ideals in rings
are assumed two-sided. 
The {\em rank} of a
ring $R$, denoted by $\dim_{\mathbb{Z}} R$,   
is the rank of its additive group, that is
$\dim_{\mathbb{Q}} R \otimes_{\mathbb{Z}} \mathbb{Q}$.

\smallskip

An algebraic closure a finite field with $q$ elements $\mathbb{F}_q$ is denoted by $\bar{\mathbb{F}}_q$.

\smallskip

We denote by $M_n(R)$ the ring of $n$-by-$n$ matrices with entries in $R$.
The subscripts in matrices and in their entries 
will always
be regarded modulo $n$.
Let $A,B \in M_n(R)$. We define
$R \langle A,B \rangle$ to be the $R$-subalgebra of $M_n(R)$ 
generated by $A$ and $B$. We will study the collection
of such pairs $(A,B)$, i.e.~the set
\[
G_n(R) = \left\{ (A,B) \in M_n(R)^2 \mid  R \langle A,B \rangle = M_n(R) \right\}.
\]

\smallskip
 
We also need the free noncommutative associative ring
$R\{ x,y \}$ whose elements we refer to as noncommutative polynomials.
The ring presentations
studied in this paper are quotients of $\mathbb{Z} \{x,y\}$. 
We do not postulate
that the identity is in $R \langle A,B \rangle$, while
we postulate that $1 \in R\{x,y\}$.

\smallskip

Many of our considerations will be based on the following
two matrices:
\[
X =  E_{21}+E_{32}+ \ldots + E_{n,n-1}+E_{1n} ~ ~ \text{and}~ ~ 
Y = E_{11} ~\text{for} ~ ~  n \geq 2 \text{.}
\] 

\smallskip

Let $FS(x,y)$ be a free semigroup on $x$ and $y$. It has
the lexicographic order as well as the word length $l(w)$
counting the total number of $x$ and $y$ in $w \in FS(x,y)$.

\paragraph{The matrices $T_{m,n,R}$.}

\smallskip

Let $R$ be a ring, and let $x_{ij}, y_{ij}$, where
$1 \leq i,j \leq n$,
be algebraically independent transcendental variables over $R$.
We see that 
$
\# \{ w \in FS(x,y) \mid l(w) \leq m\} = 
2^{m+1} -2$.
Below, we define the matrix 
\[
T_{m,n,R} \in M_{\left( 2^{m+1}-2 \right) \times n^2} \left( R[x_{ij}, y_{ij}] \right) \text{.}
\]
Let $w  = w(x,y) \in FS(x,y)$. 
We
substitute the matrices $\left( x_{ij} \right)$ and $\left( y_{ij} \right)$
for $x$ and $y$, respectively. The result is the $n$-by-$n$ matrix
$
\left( z_{ij} \right) = w_R\left( \left( x_{ij} \right), \left( y_{ij} \right) \right)
$, 
which we write 
as a row vector as follows
\begin{equation}\label{Flat}
(z_{11},z_{12} \ldots, z_{1n}, z_{21},z_{22} \ldots, z_{2n}, \ldots, z_{n1},z_{n2} \ldots, z_{nn})
\text{.}
\end{equation}
We call the operation of transforming the matrix $\left( z_{ij} \right)$ into the vector
(\ref{Flat}) {\em flattening} of $\left( z_{ij} \right)$.
We define $T_{m,n,R}$ to be the matrix whose rows are
the flattened matrices
$
w_R\left( \left( x_{ij} \right), \left( y_{ij} \right) \right)
$ such that $l(w) \leq m$,
the words $w$ being ordered lexicographically.

If $A,B \in M_n(R)$, then $T_{m,n,R}(A,B)$ is the matrix obtained from
$T_{m,n,R}$ by substituting the entries of $A$ and $B$ for $\left( x_{ij} \right)$
and $\left( y_{ij} \right)$, respectively.

\medskip

Let $S \subseteq \mathbb{Z}^m$ and $B_k = \{ (x_1, \ldots, x_m) \in 
\mathbb{Z}^m : \max_{1 \leq i \leq n} |x_i| \leq k \}$. The {\em asymptotic density}
of $S$ in $\mathbb{Z}^m$ is 
\[
\lim_{k \to \infty} \frac{\# B_k \cap S}
{\# B_k} \text{.}
\]

\subsection{Motivation and description of the main results}
The properties of the ring $M_n(\mathbb{Z})$ are based entirely on the presentation
by the elementary matrices $E_{ij}$ subject to the relations $E_{ij}E_{kl}=\delta_{jk}E_{il}$. 
This set of $n^2$ generators may be further reduced. 
Moreover, the matrices $X$ and $Y$
generate $M_n(\mathbb{Z})$. These matrices will be used
to construct several presentations of $M_n(\mathbb{Z})$ with $2$ 
generators and finitely many relations. We investigate
the interdependence between the relations in these presentations.
We also use them to construct $2$-generator presentations
with finitely many relations of certain direct sums of matrix rings.
Burnside's Theorem from \cite{Burnside1}
implies that the set $G_n(\mathbb{C})$ is infinite.
This paper, in contrast, studies
the set $G_n(\mathbb{Z})$. In particular, we describe $G_2(\mathbb{Z})$ in the following

\medskip

\textbf{Theorem \ref{2detailed}}
{\em
Let $A, B \in M_2(\mathbb{Z})$. Then $(A,B) \in G_2(\mathbb{Z})$ if and only if the following
two conditions are satisfied:
\begin{enumerate}
\item $\gcd \left( \det A,\, \det B ,\, \det(A+B)  \right) = 1$.
\item The matrices $I_2,  A,   B,   AB$ generate $M_2(\mathbb{Z})$ as a ring.
\end{enumerate}

If $I_2,A,B$ generate $M_2\left( \mathbb{Z}\right)$ as a ring,
then their $\mathbb{Z}$-linear combinations produce $I_2,A_1,B_1$
also generating $M_2\left( \mathbb{Z}\right) $ such that
\begin{equation*}
A_1 =\left( 
\begin{array}{cc}
c & 1 \\ 
1 & 0
\end{array}
\right) ~ and ~ 
B_1 =\left( 
\begin{array}{cc}
a & 0 \\ 
b & 0
\end{array}
\right) 
\end{equation*}
where $\gcd (a, b)=1$. Moreover, the matrices $I,A_1,B_1$
generate $M_2\left( \mathbb{Z}\right) $ if and only if 
\begin{equation}\label{WiseSaid1}
a^{2}-abc - b^{2}  = \pm 1 
\text{.}
\end{equation}
The set of solutions of these equations is infinite, and when $abc \neq 0$, 
this set is effectively described
in terms of the unit group of the field $\mathbb{Q} \left( \sqrt{c^2 + 4} \right)$.
}

\medskip

We show that 
$M_n(F)^2 - G_n(F)$ is
``small" for many fields.
Namely, if $F$ is a normed field having a sequence of nonzero elements whose norms tend to zero,
then the set $G_n(F)$ is dense
in $M_n(F)^2$. We also prove that 
\[\lim_{q \to \infty} \frac{\#G_n(\mathbb{F}_q)}{\#M_n(\mathbb{F}_q)^2} = 1.\]
In contrast, the set 
$M_2(\mathbb{Z})^2 - G_2(\mathbb{Z})$ is not algebraic, and
$G_2(\mathbb{Z})$ has zero asymptotic density in
$M_2(\mathbb{Z})^2$.

\medskip

The problem of minimality of presentations in ring theory admits a number of
interpretations. For example, one may search for a presentation with
the smallest number of both generators and relations. Unfortunately, 
no technique is available to solve this problem in general. More modestly, 
one may ask whether the removal of any of the relations in a given presentation
changes the ring. We study this question
and in many cases obtain information about the structure of the 
resulting over-rings.

\smallskip

We use the following noncommutative polynomials: 
\[
r_{1,n}  =  r_{1,n}(x) = x^n - 1, ~ \,
r_{2,n}  =  r_{2,n}(x,y) = \sum_{i=0}^{n-1} x^{n-i}yx^i - 1, \]
\[
s_0  =  s_0(y) = y^2 - y, ~ \,
s_j   =  s_j(x,y) = y x^j y ~ \, \text {for} ~ \, j \geq 1 \text{.}
\]  
Here are the presentations studied in our paper:
\begin{equation}\label{dvoika}
M_2(\mathbb{Z}) \cong \langle x,y \mid x^2 = y + xyx  = 1, \, yxy = 0 \rangle.
\end{equation}
\begin{equation}\label{troika}
M_3(\mathbb{Z}) \cong \langle x,y \mid x^3 = y + x^2yx + xyx^2  = 1, \, yxy = 0 \rangle.
\end{equation}
\begin{equation}\label{chetverka}
M_4(\mathbb{Z}) \cong 
\langle x,y \mid r_{1,4} = r_{2,4}  = s_0 = s_1 = 0 \rangle.
\end{equation}
\begin{equation}\label{pyaterka}
M_5(\mathbb{Z}) \cong 
\langle x,y \mid r_{1,5} = r_{2,5}  = s_0 = s_1 = 0 \rangle.
\end{equation}
\begin{equation}\label{NovyyGod}
M_n(\mathbb{Z}) \cong \langle x,y \mid r_{1,n} = r_{2,n}  = s_j = 0, ~ 
 1 \leq j \leq n-1  \rangle.
\end{equation}
\begin{equation}\label{dTolya}
M_n(\mathbb{Z}) \cong \langle x,y \mid r_{1,n} = r_{2,n}  = s_0 = s_k = 0, 
~ 1 \leq k \leq \left\lfloor n/2 \right\rfloor \rangle.
\end{equation}

While we cannot completely answer the question of
minimality in the presentations above, some information
is available in Theorems \ref{Elim}, \ref{SaidWants}, and \ref{Said3} below. 
Theorems \ref{Elim} and \ref{SaidWants} investigate the effect of the removal of 
certain relations
from (\ref{NovyyGod}).

\smallskip

\textbf{Theorem \ref{SaidWants}.}
{\em
\begin{enumerate}
\item The ring 
$
\mathcal{R} =  \langle x, y \mid r_{1,n} = s_m = 0,~
0 \leq m \leq n-1  \rangle
$
is isomorphic to a direct
sum of the rings $M_n(\mathbb{Z})$ and $\mathbb{Z}[x]/(x^n-1)$.

\item Let $\emptyset \neq H \subsetneqq N=\left\{ 1,2,...,n-1\right\} $ and
$H^{\prime }=N-H$. Suppose that $H$ satisfies the following conditions modulo $n$:
\begin{enumerate}
\item\label{good1} $\{ a+b \mid a,b \in -H \cup H \}\subseteq H^{\prime }$. 
\item\label{good2} If $ h,\,k,\,l,\, -h+k+l\in H$, then $h=k$ or $h=l$. 
\end{enumerate}
Then the ring 
$S(H) = \langle x,y \mid r_{1,n} = r_{2,n} = s_j = 0,~  j\in H^{\prime} \rangle$
has finite rank.
\end{enumerate}
}

\smallskip

\textbf{Theorem \ref{Said3}}
{\em
The ring $\mathbb{Z} \{ x,y\}$ has a quotient $R = R_n$ such that
\begin{enumerate}
\item $R$ is an over-ring of $M_n(\mathbb{Z})$.
\item Under the natural epimorphism $R \twoheadrightarrow M_n(\mathbb{Z})$, 
the images of the ideals generated by $r_{1n}, s_1, \ldots, s_n$ form a direct sum.
\end{enumerate}  
}

\smallskip

\noindent{In the proof of this theorem we introduce  an analog of 
the Magnus Embedding (see Lemma on p.~764 of Magnus \cite{Magnus}).}

\medskip

We prove the following theorem about linear representations of matrix rings.

\smallskip

\textbf{Theorem \ref{Said1}.}
{\em
Let $\mathcal{D}$ be a commutative domain of characteristic either zero or 
at least $m+1$,
over which every finitely generated 
projective module is free.
Let
$\mathcal{S} $ 
be a subring of $M_{m}\left( \mathcal{D}\right)$ generated 
by some nonzero $X_1$ and $ Y_1 $
such that
\begin{equation*}
X_1^{n+1} = X_1, ~ Y_1 X_1^n = Y_1, ~ Y_1^2 = Y_1,~
\sum_{i = 0}^{n-1} X_1^{n-i} Y_1 X_1^{i} = X_1^{n}
\text{.}
\end{equation*}
Then the trace $k$ of $Y_1$ is a positive integer, and
there exist $B \in GL_m(\mathcal{D})$ such that, putting
$r = m - kn$, we have
\begin{equation*}
B^{-1}X_{1}B=\left( 
\begin{array}{cc}
I_{k}\otimes X & 0_{k\times r} \\ 
0_{r\times k} & 0_{r\times r}
\end{array}
\right)~\, \text{and}~ \,
B^{-1}Y_{1}B=\left( 
\begin{array}{cc}
I_{k}\otimes Y & 0_{k\times r} \\ 
0_{r\times k} & 0_{r\times r}
\end{array}
\right) \text{.}
\end{equation*}
}

\smallskip

The rigidity of the embeddings in
the above theorem also follows from more general results in Azumaya algebras
(see Faith \cite{Faith}, pp.~481-482).

\medskip

We investigate the matrices satisfying the relations
of ($\ref{dTolya}$).
Let $x_1, \ldots, x_n$ be numbers. These numbers determine
the circulant matrix 
$\text{circ} (x_1, \ldots, x_n) = \sum_{i=1}^{n} x_{n-i+1} X^i$.
Integral $n$-by-$n$ circulant matrices are exactly the
elements of the group ring $\mathbb{Z} \langle X \rangle$.

\smallskip

\textbf{Theorem \ref{Said-Bogdan}.}
{\em The set 
$\mathcal{Y} = \{Y_1 \in M_n(\mathbb{Z}) \mid  Y_1^2 = Y_1, ~ r_{2,n}(X,Y_1) = 0\}$ 
has the property that the pair $(X, Y_1)$  
satisfies all relations of (\ref{dTolya}) and 
all $Y_1$ have trace $1$.
If $n = 2,3,4,6$ then 
$Y_{1} = E_{ii}$ for some $i$.
Otherwise, $\mathcal{Y}$
is
infinite, and 
if $Y_1 \neq E_{ii}$
then it has both positive and negative entries.

Any $Y_1$ is of the form $\left(c_i d_j \right)$ for some
integers $c_i, d_j$ such that
the matrices $\text{circ} \left(c_1, \ldots, c_n  \right)$ and 
$\text{circ} \left(d_1, \ldots, d_n  \right)$ are mutually inverse.
Any $Y_1$ is conjugate to $Y$ by an integral circulant matrix
with determinant $\pm 1$.
}

\smallskip

\noindent{This result depends on a classic theorem of G.~Higman \cite{Higman} about the structure
of the unit group of an integral group ring of a finite Abelian group.

\medskip

In the final part of this paper, we obtain some
$2$-generator presentations with finitely
many relations for finite direct sums of $M_n(\mathbb{Q})$, and
for the direct sums 
$\bigoplus_{j=1}^{k} M_{n_j} (\mathbb{Z})$ where
$n_1, \ldots, n_k \geq 2$ are pairwise relatively prime.

\medskip

\textbf{Acknowledgments.} 
The first author is grateful to Rostislav I.~Grigorchuk for asking very interesting questions 
leading 
to this research.
The first author also thanks Everett C.~Dade, Ronald G.~Douglas, Leonid Fukshansky, 
Gerald J.~Janusz, Doug Hensley,
Matthew Papanikolas,
Derek J.S.~Robinson,
David J.~Saltman, and Jeffrey D.~Vaaler for very useful comments and discussions.

\section{On the structure of $G_n(\mathbb{Z})$}

The starting point of this paper is the following theorem
of W.~Burnside (Burnside \cite{Burnside1}). 
We state it in the modern form, similar to Lam \cite{Lam}, p.~103.

\begin{theorem}[Burnside's Theorem]\label{Burnside's} Let $F$ be a field, $V$ a finite-dimensional 
$F$-linear space, and $S$ an $F$-subalgebra
of the algebra $End_F V$ of linear operators. Suppose that $V$ is a simple left $S$-module such that
$End_S V$ consists exactly of scalar multiples
of the identity operator on $V$.
Then $S = End_F V$.
\end{theorem}

The condition $End_S V = F \, id_{V}$ may not always be omitted if $F$ is not algebraically
closed -- counter-examples exist for any such a field. 
If $F$ is algebraically closed, however, this condition is superfluous
by Schur's Lemma (see Curtis and Reiner \cite{Curtis}, 27.3).
Burnside has proved his result in a different form from first principles by linear algebra:
see Burnside \cite{Burnside1}, p.~433, {\em Theorem}.  

In this paper, Burnside's Theorem is applied to $2$-generator subalgebras
of of $End_F V$. Therefore, below we restate the theorem for this case.

\begin{theorem}\label{Main} $F \langle A,B \rangle = End_F V$ if and only if the 
following conditions are satisfied:
\begin{enumerate}
\item\label{c1} The only subspaces of $V$, invariant under both $A$ and $B$, are $0$ and $V$.
\item\label{c2} Only scalar multiples of $id _V$ commute with both $A$ and $B$.
\end{enumerate}
\end{theorem}

\medskip

We need the following lemma that sometimes
makes it unnecessary to verify
Condition \ref{c2} of Theorem \ref{Main}. 

\begin{lemma}\label{Main1}
Let $L/F$ be a field extension, then 
$G_n(L) \cap M_n(F)^2 = G_n(F)$.
\end{lemma}

\begin{proof}
1. The inclusion $G_n(L) \cap M_n(F)^2 \subseteq G_n(F)$ holds
because linear independence over $L$ implies linear independence over $F$.

2. Conversely, let $(A,B) \in G_n(F)$. Then there exist
$n^2$ words $w_1, \ldots, w_{n^2}$ in $A,B$ that form an $F$-basis
of $M_n(F)$. It follows that $w_1, \ldots, w_{n^2}$ form
an $L$-basis of $M_n(L)$. Indeed, $E_{ij}$ form an $L$-basis of $M_n(L)$,
and the two bases are related by an invertible matrix with entries in 
$F \subseteq L$.
\end{proof}

David Saltman \cite{Saltman} has kindly communicated to us  
the following local-global principle. To state it, we need the map
$\widehat{p}: M_n(\mathbb{Z}) \to M_n(\mathbb{F}_p)$ that
reduces modulo $p$ every entry of a matrix. 

\begin{theorem}\label{Local-Global}
$G_n(\mathbb{Z}) = \bigcap_{p \, \text{prime}} \widehat{p}\,^{-1} 
\left( G_n \left( \mathbb{F}_p \right) \right)$.
\end{theorem}

\begin{proof}
We regard $M  = M_n(\mathbb{Z})$ as an additive Abelian group of rank $n^2$.
Consider the subgroup $G = \mathbb{Z} \langle A,B  \rangle$.
If $G$
is generated by $t$ elements, then their $\widehat{p}$-images
generate $\widehat{p} G$, so that
$t \geq \dim_{\mathbb{F}_p} \widehat{p}G = n^2$.
Therefore $t = n^2$, so that the index
$k = |M : G|$ is finite.

It remains to see that $k = 1$. Suppose that $k \geq 2$. We may choose
a subgroup $H$ of $M$ such that $G \subseteq H$
and $h = |M : H|$ is prime. Then $h M \subseteq H$. 
Therefore 
$|\mathbb{F}_{h}^{n^2}:\widehat{h}H| =  
|M / h M : H / h M|  
= |M : H| = h$, 
so that
$\mathbb{F}_h^{n^2} = \widehat{h} G \subseteq \widehat{h} H
\subsetneqq \mathbb{F}_{h}^{n^2}$, a contradiction.
\end{proof}

\smallskip

Combining Schur's Lemma, Lemma \ref{Main1}, Theorems \ref{Main} and \ref{Local-Global}
provides a simple method of constructing
infinitely many elements $(A,B)$ in $G_n(\mathbb{Z})$
without finding the corresponding $f_{ij} \in \mathbb{Z} \{ x,y \}$
such that $E_{ij} = f_{ij}(A,B)$.

\begin{theorem}\label{Local-Global1}
$(A,B) \in G_n(\mathbb{Z})$ if and only if  
$\bar{\mathbb{F}}_p\langle \widehat{p}A, \widehat{p}B \rangle x  = \bar{\mathbb{F}}_p^n$
for any $0 \neq x \in \bar{\mathbb{F}}_p^n$ and any prime $p$.
\end{theorem}

\begin{example}\label{Beauty1} $(X, E_{st}) \in G_n(\mathbb{Z})$ for any $s$ and $t$.
\end{example}

\begin{proof}
We apply Theorem \ref{Local-Global1}. Let 
$ x = \left(\alpha_1, \ldots, \alpha_n \right) = \sum_{i=1}^n \alpha_i e_i
\in \bar{\mathbb{F}}_p^n$ 
be a nonzero column vector. By several
applications of $X$ to $x$, we may assume that $\alpha_t \neq 0$. 
Then $y  = \alpha_t^{-1}Yx = e_s$ and
$\{X^iy \mid 1 \leq i \leq n \} = \{ e_1, \ldots, e_n \}$.
\end{proof}

\begin{example}\label{Beauty3}
Let $A = (a_{ij}),~B = (b_{ij}) \in M_n(\mathbb{Z})$ be such that 

1. $a_{l-1,l} = 1$ for $2 \leq l \leq n $
and $a_{ij} = 0$ if $i \leq j$.

2. $\{ e_1 \} \cup \{B^l e_1 \mid  2 \leq l \leq n \}$
form a $\mathbb{Z}$-basis of $\mathbb{Z}^n$.

Then $(A,B) \in G_n(\mathbb{Z})$.  
\end{example}

\begin{proof}
Let 
$ x 
\in \bar{\mathbb{F}}_p^n$ 
be nonzero, and
$k$ be the largest subscript corresponding to a nonzero component of $x$.

Case 1. If $k = 1$, then $e_1 \in \bar{\mathbb{F}}_p \langle A,B \rangle x $, so that
$\{ e_1 \} \cup \{B^l e_1 \mid  2 \leq l \leq n \}$
form a $\bar{\mathbb{F}}_p$-basis of $\bar{\mathbb{F}}_p^n$.

Case 2. If $k \geq 2$, then 
$ A^{k-1}x $
has the property that its first component is nonzero and all others are zero,
so that we return to Case 1.
\end{proof}

These examples clearly imply that the set $G_n(\mathbb{Z})$ is infinite.
This also follows from the fact that the set
$\{ (U^{-1}X U, U^{-1}Y U) \mid U \in GL_n(\mathbb{Z}) \}$ 
is infinite. Indeed,
the centralizers of $X$ and $Y$ have the following properties:
$C_{M_n(\mathbb{Z})}(X) = \mathbb{Z} \langle X  \rangle\,$, and
$C_{M_n(\mathbb{Z})}(Y)$ consists of the matrices
$(a_{ij})$ such that $a_{j1} = a_{1j} = 0$ for all $2 \leq j \leq n$.
Therefore, the intersection the two centralizers with $GL_n(\mathbb{Z})$ is $\{ \pm I_n \}$.

\medskip

Let $R$ be a commutative ring.
Following Longstaff \cite{Long}, we introduce
{\em the minimum spanning length} $\text{msl}_R$ for every $(A,B) \in G_n(R)$. 
Namely, if $(A,B) \in G_n(R)$, then
$\text{msl}_{R}(A,B)$ is the smallest integer $s$
with the property that there exist $w_1, \ldots, w_{n^2} \in FS(x,y)$ with
$\max_{1 \leq j \leq n^2} l(w_j) \leq s\,$, such that
$M_n(R) = w_{1}(A,B)R + \ldots + w_{n^2}(A,B)R$.
In the case
of fields, Proposition 1 of Longstaff \cite{Long}
is easily generalized to

\begin{lemma}\label{Hilbert-Style}
Let $F$ be a field. Then
\begin{equation}\label{EssLongstaff}
\max_{(A,B)\in G_n(F)} \text{msl}_F (A,B) \leq n^2 - 1.
\end{equation}
\end{lemma} 

\begin{proof}
Let
$\mathcal{W}_k$ be the $F$-linear span of all matrices that may be written
as $A,B$-words of length $\leq k$. We see that 
$\mathcal{W}_k \subseteq \mathcal{W}_{k+1}$. Let $m$ be the smallest value of the subscript
stabilizing this chain.
Then $\dim_F \mathcal{W}_1  = 2$, and 
$\dim_F \mathcal{W}_{l+1}  - \dim_F  \mathcal{W}_{l} \geq 1$ for any $l  \leq m-1$.
Therefore $m \leq n^2 - 1$.
\end{proof}

We extend this result to $\mathbb{Z}$ below.

\begin{theorem}\label{Revelation}
Let $A,B \in M_n(\mathbb{Z})$. Then $(A,B) \in G_n(\mathbb{Z})$ if
and only if the rows of the matrix $T_{n^2-1, n^2,\mathbb{Z}}(A,B)$
span $M_n(\mathbb{Z})$.
\end{theorem}

\begin{proof} It suffices to prove that the condition is necessary.
Let $(A,B) \in G_n(\mathbb{Z})$. Then $(A,B) \in G_n(\mathbb{F}_p)$ for 
every prime $p$. Therefore
by Lemma \ref{Hilbert-Style},
there exists a nonzero $n^2$-by-$n^2$ minor of 
$T_{n^2-1,\, n^2,\,\mathbb{F}_p}(\widehat{p} A,\, \widehat{p} B)$.
Let $w_1, \ldots w_{n^2} \in FS(x,y)$ be the words giving rise to this minor, and
let $H_p = \sum_{k=1}^{n^2} w_{k}(A,B) \mathbb{Z}$. Then
the group $H = \sum_{p~\text{prime}} H_p$ has the property
that $\widehat{p} H = M_n(\mathbb{F}_p)$ for every prime $p$.
At the same time, $H$ is a subgroup of the group generated by all row-vectors
of $T_{n^2-1,\, n^2,\,\mathbb{Z}}(A,B)$. It remains to apply
Theorem \ref{Local-Global} and Lemma \ref{Hilbert-Style}.
\end{proof}

The inequality (\ref{EssLongstaff}) is not sharp, even for $n = 2$, because 
Proposition 2 on p.~250 of Longstaff \cite{Long} implies
$
\max_{_{(A,B)\in G_2(\mathbb{C})}} \text{msl}_{\mathbb{C}}(A,B) = 2$.
This is true over any field: to modify the proof, in the last paragraph of Lemma 1 
of Longstaff \cite{Long}, we propose to replace taking adjoints with taking transposes.
The paper Longstaff \cite{Long} contains an intriguing and well substantiated conjecture 
that $\max_{_{(A,B)\in G_n(\mathbb{C})}} \text{msl}_{\mathbb{C}}(A,B) \leq 2n-2$.

\subsection{Description of $G_2(\mathbb{Z})$}

We relate below
the elements of $G_2(\mathbb{Z})$ to the solutions of the Diophantine equation
(\ref{WiseSaid1}).

\begin{theorem}\label{2detailed}
Let $I = I_2$ and $A, B \in M_2(\mathbb{Z})$. 
Then $(A,B) \in G_2(\mathbb{Z})$ if and only if the following
two conditions are satisfied:
\begin{enumerate}
\item\label{cond1} $\gcd \left( \det A,\, \det B ,\, \det(A+B)  \right) = 1$.
\item\label{cond2} The matrices $I,  A,   B,   AB$ generate $M_2(\mathbb{Z})$ as a ring.
\end{enumerate}
If $I,A,B$ generate $M_2\left( \mathbb{Z}\right)$ as a ring,
then their $\mathbb{Z}$-linear combinations produce $I,A_1,B_1$
also generating $M_2\left( \mathbb{Z}\right) $ such that
\begin{equation*}
A_1 =\left( 
\begin{array}{cc}
c & 1 \\ 
1 & 0
\end{array}
\right) ~ and ~ 
B_1 =\left( 
\begin{array}{cc}
a & 0 \\ 
b & 0
\end{array}
\right) 
\end{equation*}
where $\gcd (a, b)=1$. Moreover, the matrices $I,A_1,B_1$
generate $M_2\left( \mathbb{Z}\right) $ if and only if 
\begin{equation}\label{WiseSaid1}
a^{2}-abc - b^{2}  = \pm 1 
\text{.}
\end{equation}
The set of solutions of these equations is infinite, and when $abc \neq 0$, 
this set is effectively described
in terms of the unit group of the field $\mathbb{Q} \left( \sqrt{c^2 + 4} \right)$.
\end{theorem}

\begin{proof} Put
$\mathcal{S} = \mathbb{Z} \langle A,B \rangle$.
The Cayley-Hamilton Theorem
successively applied to the matrices $A,B,A+B$ yields 
$
\det(A)I,\, \det(B)I,\, \det(A+B) I \in \mathcal{S}$.
Since in addition, $BA = (A+B)^2 - A^2 - B^2 - AB =$ \\
\[\text{tr}(A+B)(A+B) - \det(A+B)I - 
\text{tr}(A)A + \det (A) I - \text{tr}(B)B + \det (B) I - AB\, ,\] 
we conclude that
\begin{equation}\label{Alg}
\mathcal{S} = g \mathbb{Z}I +
\mathbb{Z} A + \mathbb{Z} B + \mathbb{Z} AB,~\text{where}~
g = \gcd \left( \det A,\, \det B ,\, \det(A+B)  \right) \text{.}
\end{equation}
If $g \geq 2$, then reducing (\ref{Alg}) modulo $g$, we obtain a contradiction 
for reasons of cardinality. Therefore $g = 1$, and
$\mathcal{S} = M_2(\mathbb{Z})$ if and only if Conditions \ref{cond1} and \ref{cond2} above
are satisfied.

\medskip

Now suppose that $I,A,B$ generate the ring $M_2(\mathbb{Z})$. Let
\begin{equation*}
A=(x_{ij}),~B=(y_{ij}).
\end{equation*}
Since $ I,A,B $ generate $M_2\left( \mathbb{Z}\right) $
modulo any integer $m$, we conclude that 
$
\gcd \left( x_{12},y_{12}\right)  =1.
$ 
Let $a,b$ be integers such
that $ax_{12}+by_{12}=1$. Then 
\[
A^{\prime } = aA+bB=\left( 
\begin{array}{cc}
x_{11}^{\prime } & 1 \\ 
x_{21}^{\prime } & x_{22}^{\prime }
\end{array}
\right) ,~
B^{\prime } = B-y_{12}A^{\prime }=\left( 
\begin{array}{cc}
y_{11}^{\prime } & 0 \\ 
y_{21}^{\prime } & y_{22}^{\prime }
\end{array}
\right)
\]
and therefore $ I,A^{\prime },B^{\prime } $ generate 
$M_2\left(\mathbb{Z}\right) $. We use the identity $I$ to obtain
\[
A^{\prime \prime } = A^{\prime }-x_{22}^{\prime }I=\left( 
\begin{array}{cc}
x_{11}^{\prime \prime } & 1 \\ 
x_{21}^{\prime } & 0
\end{array}
\right) , ~
B^{\prime \prime }  = B^{\prime }-y_{22}^{\prime }I=-\left( 
\begin{array}{cc}
y_{11}^{\prime \prime } & 1 \\ 
y_{21}^{\prime } & 0
\end{array}
\right).
\]
Again, $ I,A^{\prime \prime },B^{\prime \prime } $ generate 
$M_2\left( \mathbb{Z}\right) $. We rewrite 
$A^{\prime \prime },B^{\prime \prime }$ as $A,B$, respectively; that is, we may assume 
\begin{equation*}
A=\left( 
\begin{array}{cc}
x_{11} & 1 \\ 
x_{21} & 0
\end{array}
\right) , ~B=\left( 
\begin{array}{cc}
y_{11} & 0 \\ 
y_{21} & 0
\end{array}
\right) \text{.}
\end{equation*}
Let $c,d$ be integers such that $cx_{21}+dy_{21}=1$. We may replace $A$ by 
\[
A^{\prime }=cA+dB=\left( 
\begin{array}{cc}
x_{11}^{\prime } & 1 \\ 
1 & 0
\end{array}
\right). 
\]
Thus we may assume that
\begin{equation*}
A=\left( 
\begin{array}{cc}
x_{11} & 1 \\ 
1 & 0
\end{array}
\right) , ~B=\left( 
\begin{array}{cc}
y_{11} & 0 \\ 
y_{21} & 0
\end{array}
\right) \text{,} ~ \gcd \left( y_{11},y_{21}\right) =1.
\end{equation*}

We want to determine when the $\mathbb{Z}$-span of $I$,
$A$, $B$, $AB$
is $M_2(\mathbb{Z})$. If $E_{11}$ is a linear combination of $I,A,B,AB$ 
then $E_{12} + E_{21} \in \langle I,A,B \rangle\,$, and 
therefore $\langle I,A,B \rangle = M_2\left(\mathbb{Z}\right) $.

Let $a,b,c,d$ be integers such that $aI+bA+cB+dAB= E_{11}$. As 
\begin{equation*}
aI+bA+cB+dAB=\left( 
\begin{array}{cc}
a+bx_{11}+cy_{11}+d\left( x_{11}y_{11}+y_{21}\right)  & b \\ 
b+cy_{21}+dy_{11} & a
\end{array}
\right) \, \text{,}
\end{equation*}
the above equation has a solution if and only if 
\begin{equation*}
a=b=0,~dy_{11}=-cy_{21}, ~ cy_{11}+d\left( x_{11}y_{11}+y_{21}\right) =1.
\end{equation*}
If $y_{11}=0\,$,  then $dy_{21}=1\,$; therefore $y_{21}=d=\pm 1$ and $c=0.$
Similarly if $y_{21}=0\,$, then $y_{11}=c=\pm 1$ and $d=0$.

We assume $y_{11}, y_{21} \neq 0$. Therefore $c, d \neq 0\,$,  and since 
$\gcd \left( y_{11},y_{21}\right) =1$, from $dy_{11}=-cy_{21}$ we conclude
that there exists an integer $c^{\prime }$ such that
\begin{equation*}
c=c^{\prime }y_{11},~d=-c^{\prime }y_{21}\text{.}
\end{equation*}
The equation $cy_{11}+d\left( x_{11}y_{11}+y_{21}\right) =1$ yields
$c^{\prime }\left( y_{11}^{2}-x_{11}y_{21}y_{11}-y_{21}^{2}\right) =1$
therefore
$y_{11}^{2} - x_{11}y_{21} y_{11} - y_{21}^{2} = \pm 1$.
It remains to write $a = y_{11}$, $b = y_{21}$, $c = x_{11}$, and we obtain
(\ref{WiseSaid1}).
Since it is easy to analyze the solutions
when one of $a,b,c$ is zero, we will investigate
the other solutions only. Equation (\ref{WiseSaid1}) is quadratic in $a$;
therefore, a necessary condition for (\ref{WiseSaid1}) to have integral
solutions is that the equation 
\begin{equation}\label{Disc}
d^2 = (bc)^2 + 4 (b^2 \pm 1)
\end{equation}
should have integral solutions too. If this is so, then
\begin{equation}\label{NaiveQuad} 
a = \frac{bc \pm d}{2}
\end{equation}
From (\ref{Disc}) we observe that $d \equiv d^2 \equiv (bc)^2 \equiv bc \, (\text{mod} \, 2)$.
In other words, (\ref{Disc}) implies (\ref{NaiveQuad}).
Now (\ref{Disc}) may be rewritten as
\begin{equation}\label{Pell}
d^2 - (c^2 + 4) b^2 = \pm 4.
\end{equation}
Let $s$ be the square-free part of the number $c^2 + 4$. Then
according to Fr\"{o}hlich and Taylor \cite{Frohlich}, 1.3, the units of 
$\mathbb{Q} \left( \sqrt{c^2 + 4}  \right)$
uniquely, under the map $(d,b) \mapsto (1/2) (d + b \sqrt{c^2+4})$, correspond to the integral solutions
of (\ref{Pell}). There are infinitely many of them by the Dirichlet's Unit Theorem. 
Algorithm 5.7.2 in Cohen \cite{Cohen} computes the fundamental unit of a rational
quadratic number field with positive discriminant. 

Therefore, for a fixed $c$, we can produce units in $\mathbb{Q} \left( \sqrt{c^2 + 4}  \right)$,
thus determining $b$ and $d$; then $a$ may found from (\ref{NaiveQuad}). 
\end{proof}

\subsection{Asymptotic properties of $G_n(\mathbb{Z})$}

\begin{lemma}\label{Diophantine1} 
Let 
$0 \neq f \in \mathbb{Z}[x_1, \ldots, x_n]$. Then
$V(f) = \{ a \in \mathbb{Z}^n \mid f(a) = 0 \}$ has zero asymptotic density
in $\mathbb{Z}^n$.
\end{lemma}

\begin{proof} 
Put $B_k = \{ (a_1, \ldots, a_n) \in \mathbb{Z}^n \mid -k \leq a_i \leq k~\text{for}~\text{all}~i \}$.
The case $n = 1$ is clear because
$\# B_k \leq \deg(f)$ for all $k$.

Let $n = 2$, $x = x_1$, $y = x_2$ and $d = \deg (f)$. Then
$f(x,y) = \sum_{j=1}^d f_j(x) y^j$ for some $f_j(x) \in \mathbb{Z}[x]$.
Let 
$
S = \{-k \leq a \leq k \mid f_j(a) = 0~\text{for all}~j \}
$.
Then $\# S \leq d$. We may write
$V(f) = A \cup B$,
where
\[A = \{(a,b) \in V(f) \mid a \in S \}~\,
\text{and}~\, B = \{(a,b) \in V(f) \mid a \notin S \} \text{.}
\]
If $a \in \{-k, \ldots, k \} - S$,
then $\# \{b \mid  (a,b) \in B \} \leq d $.
Hence, 
\[
\# V(f) \leq \# A + \# B \leq
\# S \, \# \{-k ,\ldots, k\} + \# \left( \{-k ,\ldots, k\} - S \right) \, d = O(k).
\]
Since $\# B_k = (2k+1)^2$, we conclude that the lemma is true when $n = 2$.

The case $n \geq 3$ is handled similarly by induction on $n$.
\end{proof}

The exponent of $k$ in the estimate 
$  \# \left( B_k \cap V(f) \right)/
 \# B_k   = O \left(k^{-1} \right)$
in the proof of Lemma \ref{Diophantine1} is the best possible in general, as exemplified
by the polynomial $f(x_1, \ldots, x_n) = x_1$.

\begin{corollary}\label{Diophantine}
The set $M_n(\mathbb{Z})^n - G_n(\mathbb{Z})$ is not algebraic.
\end{corollary}

\begin{proof}
Suppose that the theorem is false.
Then Lemma \ref{Diophantine1}
implies that $G_n(\mathbb{Z})$ has asymptotic density $1$ in $M_n(\mathbb{Z})^2$.
This is false, however, because 
$M_n(2\mathbb{Z})^2 \subseteq M_n(\mathbb{Z})^2 - G_n(\mathbb{Z})$,
and $M_n(2\mathbb{Z})^2$ has  asymptotic density $2^{-2 n^2}$ in $M_n(\mathbb{Z})^2$, implying that
$G_n(\mathbb{Z}) \cap M_n(2 \mathbb{Z})^2$
is non-empty.
\end{proof}

\begin{theorem}\label{Soft} The set $G_2(\mathbb{Z})$ has zero asymptotic density in $M_2(\mathbb{Z})^2$.
\end{theorem}

\begin{proof} Put $I = I_2$. Let 
$A,B \in M_2(\mathbb{Z})$ such that
$I$, $A$, $B$ generate $M_2(\mathbb{Z})$ as a ring. Put 
$\mathcal{S} = \mathbb{Z} \langle A,B \rangle$.
The Cayley-Hamilton Theorem
applied to the matrices $A$, $B$, $A+B$ yields that $A^2, B^2, (A+B)^2$
are integral linear combinations of $I$, $A$, $B$. 
Since in addition, $BA = (A+B)^2 - A^2 - B^2 - AB \,$,
we conclude that
$\mathcal{S} = \mathbb{Z}I +
\mathbb{Z} A + \mathbb{Z} B + \mathbb{Z} AB$. Let
$T$ be a $4$-by-$4$ matrix whose rows are the flattened
matrices $I$, $A$, $B$, and $AB$. Then $\mathcal{S} = M_2(\mathbb{Z})$
if and only if $\det T = \pm 1$. It remains to apply
Lemma \ref{Diophantine1}.
\end{proof}

This result sometimes clarifies the relationship between $G_2(\mathbb{Z})$ and
the other subsets of $M_n(\mathbb{Z})^2$. We will give an example.
Let
$S$ be set of all
$(A,B) \in M_2(\mathbb{Z})^2 - G_2(\mathbb{Z})$
such that all the $8$ entries are relatively prime in pairs. We will see that
asymptotically, almost all elements of $S$ lie outside of $G_2(\mathbb{Z})$. To 
formalize this statement,
let
$m_k = \prod_{p \, \text{prime}, \, p \leq k} p$ and
$
D_k = \left\{ (a_1, \ldots, a_8) \in \mathbb{Z}^8 : 
\max_{1 \leq i \leq 8} |a_i| \leq m_k \right\}$.
We claim that
\begin{equation}\label{Probability}
\lim_{k \to \infty} \frac{\# S \cap D_k}{\# D_k} = \prod_{p \, \text{prime}} (p-1)^7 (p+7) p^{-8}>0.
\end{equation}
We give a heuristic argument first.
For a fixed
prime $p$, we
consider the Bernoulli scheme of choosing $8$ integers
independently and at random with the probability of success $p^{-1}$. 
Then the probability of at most $1$ success is
$\left(1- p^{-1}  \right)^8 + {8 \choose 1} p^{-1}\left(1- p^{-1}  \right)^7 = (p-1)^7 (p+7) p^{-8}$.
Taking the product over all primes gives (\ref{Probability}).

Next we prove (\ref{Probability}). We thank Doug Hensley
\cite{Hensley} for communicating the following argument to us. It is convenient
to decrease the sets $S$ and $D_k$ to retain only the $8$-tuples with all positive entries.
For a prime $p$, let $S_p$ be the set of all $8$-tuples
$(a_1, \ldots, a_8)$ whose entries are positive integers, and $p \nmid \gcd (a_i, a_j)$
if $i \neq j$. Then $S = \bigcap_{p } S_p$. The Chinese Remainder
Theorem applied to the ring $\mathbb{Z} / m_k \mathbb{Z}$ implies
\begin{equation}\label{Doug1}
\frac{\# S \cap D_k}{\# D_k} \leq 
\frac{\# \bigcap_{p \leq k} S_p \cap D_k}{\# D_k}  =   
\prod_{p \leq k} (p-1)^7 (p+7) p^{-8}.
\end{equation}
For the primes $p > k$, we have
$
\# S_p \cap D_k \leq {8 \choose 1} m_k \lfloor m_k/p \rfloor^7 $. Therefore
\begin{equation}\label{Doug2}
\frac{\# S \cap D_k}{\# D_k} \geq 
\frac{\# \bigcap_{p \leq k} S_p \cap D_k}{\# D_k} - 
\sum_{p > k} \frac{\# S_p \cap D_k}{\# D_k} = 
\prod_{p \leq k} (p-1)^7 (p+7) p^{-8} + o \left( 1 \right) \text{.}
\end{equation}
Comparing (\ref{Doug1}) and (\ref{Doug2}) yields (\ref{Probability}).

\subsection{Asymptotic and topological properties of $G_n(F)$ for fields}
\begin{lemma}\label{Density}
Let $F$ be a field. Then $M_n(F)^2 - G_n(F)$ is a non-empty algebraic set
consisting of all $(A,B) \in  M_n(F)^2$
such that the matrix $T_{n^2-1, n^2, F}(A,B)$ does not have full rank. 
\end{lemma}

\begin{proof} 
The equality of the two sets above follows
from Lemma \ref{Hilbert-Style}. The set $M_n(F)^2 - G_n(F)$
is non-empty because $G_n(\mathbb{Z})$ is non-empty.
\end{proof}

Next, we will apply Lemma \ref{Density} to normed fields satisfying the following

\begin{property}\label{Natproperty}
$F$ is a normed field (with the norm denoted by $|\cdot|$) such that
for any $\varepsilon >0$ there exists $0 \neq a_{\varepsilon} \in F$
with $|a_{\varepsilon}| < \varepsilon$.
\end{property}

Among the fields having Property \ref{Natproperty} are all the subfields
of $\mathbb{C}$ or $\mathbb{C}_p$ with their respective standard Euclidean
or $p$-adic norms.

\begin{lemma}\label{Line} 
Let $F$ have Property \ref{Natproperty}, and let 
$Z \subsetneqq F^n$ be an
algebraic set. Then $F^n - Z$ is dense in $F^n$ in the norm
topology.
\end{lemma}

\begin{proof}
Let $z \in Z$. We show that there exists
a sequence $\{ z_n\}$ in $F^n - Z$ with
$\lim_{n \to \infty} ||z - z_n|| = 0$.
Since $Z \subsetneqq F^n$, there exists a line
$L_z$ passing though $z$ and not contained in $F^n$.
Substituting the parametric equations for $L_z$ into the 
polynomial equations defining $Z$, we obtain
a system of equations in one variable, which has 
finitely many solutions, one of them being $z$.
We may choose $\varepsilon > 0$ sufficiently small
to ensure that $z$ is the only solution contained
in the ball $B_{\varepsilon}(z)$ of radius 
$\varepsilon$ and centered at $z$. Then 
there exists a sequence $\{ z_n\}$ in
$B_{\varepsilon}(z) \cap L_z$
such than $z_n \neq z$ and $\lim_{n\ \to \infty}||z-z_n|| = 0$.
In particular $z_n \in F^n - Z$.
\end{proof}

\begin{theorem}\label{Ndense}
Let $F$ have Property \ref{Natproperty}. Then
$G_n(F)$ is open and dense in $M_n(F)^2$ in the norm topology.
\end{theorem}

\begin{proof} The result follows from Lemmas \ref{Density} and \ref{Line}.
\end{proof}

Next we consider similar results for finite fields.

\begin{lemma}\label{Crutch1} Let $0 \neq f \in \mathbb{F}_q [x,y]$ and
$V(f)  = \{ v \in \mathbb{F}_q ^2 \mid f(v) = 0 \}$.
Then $ \# V(f) \leq 2 q \deg(f) $.
\end{lemma}

\begin{proof} Let $d = \deg (f)$. Then
$f(x,y) = \sum_{j=0}^d f_j(x) y^j$ for some $f_j(x) \in \mathbb{F}_q [x]$.
Let $S = \{a \in \mathbb{F}_q \mid f_0(a) = \ldots = f_d(a) = 0 \}$. Then
$\# S \leq d$. 

For every $a \in S$, there are at most $q$ values of $b \in \mathbb{F}_q$
such that $(a,b) \in V(f)$. Let
$A = \{ (a,b) \in  V(f) \mid a \in S\}$. Then $\# A \leq qd$.

Next let $B = \{ (a,b) \in  V(f) \mid a \notin S\}$.
Then there are  at most $d$ values of $b \in \mathbb{F}_q$ such that 
$(a,b) \in V(f)$ for some $a \in \mathbb{F}_q$. 
Then
$\# B \leq qd$.

Finally, $V(f) = A \cup B$, so that
$\# V(f) \leq \# A + \# B \leq 2 qd$.
\end{proof}

\begin{theorem}\label{Fdense}
For a fixed $n \geq 2$, we have
\[
\lim_{q \to \infty} \frac{\#G_n(\mathbb{F}_q)}{\#M_n(\mathbb{F}_q)^2} = 1.
\]
\end{theorem}

\begin{proof}
By Lemma \ref{Density}, $M_n(\mathbb{F}_q)^2-G_n(\mathbb{F}_q)$ is an intersection 
of finitely many hypersurfaces,
each of them having $O(q^{2n^2-1})$ points over $\mathbb{F}_q$ by Lemma \ref{Crutch1}. 
Each such a hypersurface is defined by a polynomial equation in $2 n^2$ variables with 
coefficients
in $\mathbb{Z}$, the equations being independent of $\mathbb{F}_q$.
It follows that
\[ 1 \geq \frac{\#G_n(\mathbb{F}_q)}{\#M_n(\mathbb{F}_q)^2} 
= 1 - \frac{\# \left( M_n(\mathbb{F}_q)^2 - G_n(\mathbb{F}_q) \right)}{\#M_n(\mathbb{F}_q)^2}
\geq
1 - \frac{O(q^{2n^2-1})}{q^{2n^2}} \xrightarrow[q \to \infty]~ 1.
\]
\end{proof}
However, we do not know whether the following limit exists:
\begin{equation}\label{MysteryLim}
\lim_{n,\,q \to \infty} \frac{\#G_n(\mathbb{F}_q)}{\#M_n(\mathbb{F}_q)^2} \text{.}
\end{equation}

Lemma \ref{Density} together with Theorems \ref{Main}, \ref{Ndense}, and \ref{Fdense}
imply that the set of $ (A,B) \in M_n(F)^2 $ having a proper common invariant subspace,
is small in the appropriate sense. We note that our arguments do not involve
characteristic polynomials.

\smallskip

\section{Presentations of $M_n(\mathbb{Z})$ and their applications}
We begin by recalling the definitions of the matrices $X =  \sum_{i=1}^n E_{i+1,i}$ and 
$Y = E_{11}$ for some fixed $n \geq 2$, and the noncommutative polynomials 
\[
r_{1,n}  =  r_{1,n}(x) = x^n - 1, ~ \,
r_{2,n}  =  r_{2,n}(x,y) = \sum_{i=0}^{n-1} x^{n-i}yx^i - 1,
\]
\[
s_0  =  s_0(y) = y^2 - y, ~ \,
s_j   =  s_j(x,y) = y x^j y ~ \, \text {for} ~ \, j \geq 1 \text{.}
\]

\begin{theorem}\label{Grig}
The ring $M_n(\mathbb{Z})$ has the following presentations:
\begin{equation}\label{Grigdream} 
\langle x, y \mid r_{1,n} = r_{2,n} = s_m = 0,~
1 \leq m \leq n-1  \rangle, 
\end{equation}
\begin{equation}\label{Dnepr}
\langle x, y \mid r_{1,n} = r_{2,n} = s_0 = s_k = 0,~1 \leq k \leq \lfloor n/2 \rfloor \rangle. 
\end{equation}
Both ring isomorphisms are obtained by mapping
$x$ to $X$
and $y$ to $Y$.
\end{theorem}

\begin{proof}
We see that $X$ and $Y$ satisfy all the relations of (\ref{Grigdream}) and 
(\ref{Dnepr}).

Next we prove that (\ref{Grigdream}) is a presentation
of $M_n(\mathbb{Z})$. To fix the notation, let $\mathcal{R}$ be the ring
defined by (\ref{Grigdream}).
We observe that 
\[ 
1 \cdot y = \left( \sum_{i = 0}^{n-1} x^{n-i} y x^i  \right) y =
y^2 + \sum_{i = 1}^{n-1} x^{n-i} (y x^i y)
= 
y^2
.\]
Therefore, $\mathcal{R}$ is spanned as an Abelian group
by the $n^2$ elements $x^i yx^j$ where $1 \leq i,j \leq n $;
hence $\dim_{\mathbb{Z}} \mathcal{R} \leq n^2$.
On the other hand, the map $\alpha$ given by $\alpha(x) = X$ and $\alpha(y) = Y$
extends to the ring epimorphism $\alpha: \mathcal{R} \twoheadrightarrow M_n(\mathbb{Z})$
because $E_{ij} = X^{i-1} Y Y^{1-j}$.

\smallskip

It remains to show that (\ref{Dnepr}) is a presentation
of $M_n(\mathbb{Z})$. Since all the relations of (\ref{Dnepr})
hold in (\ref{Grigdream}), it remains to establish the converse.
We propose to consider the cases of $n$ even and odd separately.
The arguments involved in either of them are the same; therefore, we will
do only the case when $n = 2s+1$ is odd. Multiplying the relation
$1 = \sum_{i=0}^{n-1} x^{n-i}yx^i$ by $y$ on the right yields 
\begin{multline}\label{riv1}
y = 1\,y = y^2+x^{n-1}(yxy)+x^{n-2}(yx^2y)+ \ldots
+ x^{n-s+1}(yx^sy)+ \\
x^{n-s}yx^{s+1}y+ \ldots+ xyx^{n-1}y.
\end{multline}
Since $y^2 = y$ and $yxy = yx^2y = \ldots = yx^sy = 0$, and $x$ is invertible, the formula
(\ref{riv1}) shortens:
\begin{equation}\label{riv2}
yx^{n-1}y + xyx^{n-2}y + \ldots + x^syx^{s+1}y = 0.
\end{equation}
Multiplying (\ref{riv2}) on the left by $y$, as before, yields
\begin{equation}\label{riv3}
yx^{n-1}y = 0,
\end{equation}
which is partly what we need.
Now substitute (\ref{riv3}) in (\ref{riv2}), cancel by $x$ on the left,
and then multiply by $y$ on the left. The result is
$yx^{n-2}y = 0$. In a similar fashion, it follows that
all $s_j(x,y) = 0$ for all $j$.
\end{proof}

The next theorem shows that Presentation \ref{Dnepr} for $n = 4,5$ may be shortened.

\begin{theorem}\label{Said45}
\begin{equation}\label{Said451}
M_4(\mathbb{Z}) \cong 
\langle x,y \mid r_{1,4} = r_{2,4} = s_0 = s_1 = 0 \rangle.
\end{equation}
\begin{equation}\label{Said452}
M_5(\mathbb{Z}) \cong 
\langle x,y \mid r_{1,5} = r_{2,5}  = s_0 = s_1 = 0 \rangle.
\end{equation}
\end{theorem}

\begin{proof}
\textbf{1.} To prove (\ref{Said451}), observe that
$0 = y r_{2,4} = s_3 x + s_2 x^2$, so that
$s_3 = -s_2 x$ and $s_3 = s_3 y  = -s_2 x y = -yx^2y(yxy) = 0$.
Therefore $s_2 = s_3 = 0$, and the result follows from Theorem \ref{Grig}.

\smallskip

\textbf{2.} We prove (\ref{Said452}) in several steps.
\[
0 = y r_{1,5}  = 
y + s_4 x + s_3 x^2 + s_2 x^3 + s_1 x^4 - y = 
s_4 x + s_3 x^2 + s_2 x^3 + s_1 x^4 =  0  \Longrightarrow 
\]
\begin{equation}\label{m51}
s_4 + s_3 x + s_2 x^2 = 0.
\end{equation} 
Similarly, by expanding $0 =  r_{2,5} y $
we have
\begin{equation}\label{m52}
s_4 + x s_3  + x^2 s_2 = 0.
\end{equation} 
Multiply (\ref{m52}) by $y$ on the right:
\begin{equation}\label{m53}
s_4 + s_2 ^2 = 0.
\end{equation}
Equate (\ref{m51}) and (\ref{m52}):
$s_3 x + s_2 x^2  = x s_3 + x^2 s_2$, and then multiply the result by $y$
on the right: $s_2^2 = x s_3 + x^2 s_2$ implying
\begin{equation}\label{m54}
s_3 = x^4 s_2 ^2  - x s_2.
\end{equation}
Multiply (\ref{m54}) by $y$ on the left 
$s_3 = y s_3 = y x^4 s_2 - yx s_2 = s_4 s_2 ^2$ and use (\ref{m53}):
\begin{equation}\label{m55}
s_3 = - s_2 ^4.
\end{equation}
Substitute (\ref{m55}) in (\ref{m55}):
\begin{equation}\label{m56}
- s_2 ^4 = x^4 s_2 ^2  - x s_2.
\end{equation}
Multiply (\ref{m56}) by $y x^2$ on the left and then use (\ref{m55}):
\[
- y x^2 s_2 ^4 = y x^2 x^4 s_2 ^2  - y x^2 x s_2 \Longrightarrow
-s_2 ^5  = - s_3 s_2  = - (- s_2 ^4) s_2 = s_2 ^5 \Longrightarrow
\]
\begin{equation}\label{m57}
2 s_2 ^5 = 0.
\end{equation}
Multiply (\ref{m56}) by $y x^4$ on the left:
$- s_4 s_2 ^4 = s_3 s_2 ^2  - s_2$. Then by
(\ref{m53}) and (\ref{m55}):
$s_2 ^6 = s_3 s_2 ^2  - s_2  = (-s_2 ^4) s_2 ^2 - s_2$.
Finally, by (\ref{m57}): $s_2 = -2 s_2 ^6  = - s_2 (2 s_2 ^5)$,
and the claim follows from Theorem \ref{Grig}. 
\end{proof}

Next we record some properties of Presentations (\ref{Grigdream}) and (\ref{Dnepr})
in connection with their minimality.
\begin{theorem}\label{Elim}
\begin{enumerate}
\item\label{it1} 
The ring
$\langle x,y \mid r_{2,n} = s_j = 0,
1 \leq j \leq n-1 \rangle$
has infinite rank.

\item\label{it2} 
The ring
$\langle x,y \mid r_{1,n} =  s_j = 0,
1 \leq j \leq n-1 \rangle$
has infinite rank. 

\item\label{it5} 
$\langle x,y \mid  r_{1,n} = r_{2,n} = 0 \rangle \ncong M_n(\mathbb{Z})$.

\item\label{it4} If $1 \leq k \leq n-1$ and $k \neq n/2$, then the relation
$s_k = 0$ follows from the other relations in (\ref{Grigdream}). 
In particular, this explains why
(\ref{troika}) is a presentation of $M_3(\mathbb{Z})$.

\item\label{it7} Removing from (\ref{Grigdream}) any two relations 
$s_h = s_{n-h} = 0$ results in a ring of an infinite rank.

\item\label{it8} Removing from (\ref{Grigdream}) any two relations 
$s_h = s_{2h} = 0$, provided $1 \leq h < 2h \leq n-1$, 
results in a ring of an infinite rank.
\end{enumerate}
\end{theorem}

\begin{proof}

\textbf{\ref{it1}.} 
Let $\mathbb{Z}(t)$ be the ring of rational functions in $t$
with integral coefficients. Consider the matrices
$A = t \sum_{i=0}^{n-1} E_{i+1,i}$ and $B = \left(1 / t^n \right) E_{11}$.
Let $\mathcal{R}$ be the subring of $M_n(\mathbb{Z}(t))$ generated
by $A$ and $B$.
These matrices satisfy all the relations of  $\mathcal{R}$. At the same time,
$A^n = t^n I \in \mathcal{R}$, so that $\mathcal{R}$ contains
$\sum_{k=1}^{\infty} t^{kn} I$, an Abelian subgroup of infinite rank.

\medskip

\textbf{\ref{it2}.} 
Consider the matrices $A = \sum_{i=0}^{n-1} E_{i+1,i}$ and 
$B = t E_{11}$. Let $\mathcal{R}$ be the subring of $M_n(\mathbb{Z}[t])$ generated
by $A$ and $B$. These matrices satisfy all the relations of 
$\mathcal{R}$. At the same time, 
$\sum_{i=0}^{n-1} A^{-i} B A^{i} = tI_n \in \mathcal{R}$, and as above, we
conclude that $\mathcal{R}$ has infinite rank.

\medskip

\textbf{\ref{it5}.} Suppose the claim is false. Then by mapping $y$ to zero, we have
$
M_n(\mathbb{Z}) \cong \langle x,y \mid r_{1,n} = r_{2,n} = 0 \rangle
\twoheadrightarrow 
\mathbb{Z}[x]/(x^n - 1)$,
but the ring $M_n(\mathbb{Z})$ does not a have proper ideal of infinite index.

\medskip

\textbf{\ref{it4}.} We need to show that the relation $y x^k y = 0$ follows from the other relations
of (\ref{Grigdream}). We have 
\[0 = r_{2,n} y =
y^2 + x^{-k}yx^ky - y~\,\text{and}~\, 
0 = y  r_{2,n}  = y^2 + y x^k y x^{-k} -y.
\]
Hence
\begin{equation}\label{m1}
y - y^2 = y x^k y x^{-k}  =  x^{-k}yx^ky.
\end{equation}

Next, we work with the expressions $y(y - y^2)$ and $(y - y^2) y$ with the help
of (\ref{m1}). 
We see that on the one hand, $y (y - y^2) =  (y x^{-k} y) x^k y = 0$, and
on the other hand $y(y - y^2) = y^2x^syx^{-k}$. Therefore $y^2x^kyx^{-k} = 0$, and since
$x$ is invertible, 
\begin{equation}\label{m2}
y^2x^ky = 0.
\end{equation}
Likewise, $(y - y^2) y = x^{-k} y x^k y^2 = y x^k(y x^{-k} y) = 0$, so that
\begin{equation}\label{m3}
yx^ky^2 = 0.
\end{equation}
Applying (\ref{m1}), (\ref{m2}), and (\ref{m3}) yields
\begin{equation*}
y x^k y = y x^k(y -  y^2) = y x^k (x^{-k} y x^k y) = y^2 x^k y = 0.
\end{equation*}

\medskip

\textbf{\ref{it7}.} It suffices to give an example of the ring of infinite rank,
where all the relations of (\ref{Grigdream}) are satisfied except for
$y x^h y = y x^{n-h} y  = 0$.

Let $\mathbb{Z}[t]$ be a polynomial ring, $X$ be the permutational matrix of
oder $n$ acting on columns, and $Y_1 = t E_{11} + (1-t) E_{1+h, 1+h}$. We denote
by $\mathcal{R}$ the subring of $M_n(\mathbb{Z}[t])$ generated by $X$ and $Y_1$. If
$1 \leq i \leq n-1$, then
\begin{equation}\label{m4}
X^i Y_1 X^{-i} = t E_{1+i, 1+i} + (1-t) E_{1+h+i, 1+h+i}
\end{equation}
implying $\sum_{i=0}^{n-1} X^i Y_1 X^{-i} = I$. 
Next, multiply (\ref{m4}) by $Y$ on the left:
\begin{equation}\label{m5}
Y_1 X^i Y_1 X^{-i} = 
t(1-t)(E_{11} E_{1+h+i,1+h+i}+ E_{1+h,1+h} E_{1+i,1+i}).
\end{equation}
We see that $Y_1 X^i Y_1 X^{-i} =0$, and therefore $Y_1 X^i Y_1 =0$, 
unless $i = \pm h$. In the latter cases we have
that $Y_1 X^h Y_1 X^{-h} = t(1-t)E_{11}$ and $Y_1 X^{-h} Y_1 X^{h} = t(1-t)E_{1+h, 1+h}$.
Therefore, in $\mathcal{R}$ all the relations of (\ref{Grigdream}) are satisfied except for
$y x^h y = y x^{n-h} y  = 0$.
Another consequence of (\ref{m5}) is $t(t-1)I_n \in \mathcal{R}$ because
$
\sum_{i=1}^{n-1} X^{-i} (Y_1 X^h Y_1 X^{-h}) X^{i} = 
t(1-t) \sum_{i=1}^{n-1} X^{-i} E_{11} X^{i} = t(t-1)I_n$.
Therefore, $\mathcal{R}$ contains an Abelian subgroup
of infinite rank, implying
that the rank of $\mathcal{R}$ is infinite as well. 

\medskip

\textbf{\ref{it8}.} As above, it suffices to give an example of the ring 
of infinite rank,
where all the relations of (\ref{Grigdream}) are satisfied except for
$yx^h = yx^{2h}y = 0$, provided $1 \leq h < 2h \leq n-1$.

Let $\mathbb{Z}[t]$ be a polynomial ring, $X$ be the permutational matrix of
order $n$ acting on columns, and $Y_1 = E_{11} + tE_{1,1+h} -t E_{1-h, 1}$.

The relation $\sum_{i=0}^{n-1} X^i Y_1 X^{-i} = I_n$ is satisfied because the subscripts
$(1,1+h)$ and $(1-h,1)$ are in the same orbit of $X$.

Next we investigate the monomial relations.
\begin{multline}\label{m6}
Y_1 X^{i} Y_1 X^{-i} = \\
\left(E_{11} + t E_{1,1+h} -t E_{1-h,1} \right)
\left( E_{1+i,1+i} + t E_{1+i, 1+h+i} - t E_{1+i-h, 1+i}\right).
\end{multline} 
On multiplying out, we see that (\ref{m6}) is zero unless $i = h,2h$. 
In the latter two cases, we have that
\begin{equation*}
Y_1 X^{h} Y_1 X^{-h} = t^2 (E_{1,1+2h} + E_{1-h, 1+h})~\text{and}~
Y_1 X^{2h} Y_1 X^{-2h} = -t^2 E_{1,1+h}.
\end{equation*}
Finally,
$ -\sum_{i=0}^{n-1} X^{-i} (Y_1 X^{2h} Y_1 X^{-2h}) X^i =
\sum_{i=0}^{n-1} X^{-i} t^2 E_{1,1+h} X^i = t^2 X^{1-h}$.
Therefore, the ring generated by $X$ and $Y_1$ has infinite rank.
\end{proof}

The above theorem describes some situations (with
the possible exception of Part \ref{it5}) where the removal 
of certain relations results in a ring of infinite rank. In contrast, the theorem
below gives two instances in which the removal of certain relations
results in a ring of finite rank.

\begin{theorem}\label{SaidWants}
\begin{enumerate}
\item The ring 
$
\mathcal{R} =  \langle x, y \mid r_{1,n} = s_m = 0,~
0 \leq m \leq n-1  \rangle
$
is isomorphic to a direct
sum of the rings $M_n(\mathbb{Z})$ and $\mathbb{Z}C_n$.

\item Let $\emptyset \neq H \subsetneqq N=\left\{ 1,2,...,n-1\right\} $ and
$H^{\prime }=N-H$. Suppose that $H$ satisfies the following conditions modulo $n$:
\begin{enumerate}
\item\label{good1} $\{ a+b \mid a,b \in -H \cup H \}\subseteq H^{\prime }$. 
\item\label{good2} If $ h,\,k,\,l,\, -h+k+l\in H$, then $h=k$ or $h=l$. 
\end{enumerate}
Then the ring 
$S(H) = \langle x,y \mid r_{1,n} = r_{2,n} = s_j = 0,~  j\in H^{\prime} \rangle$
has finite rank.
\end{enumerate}
\end{theorem}

\begin{proof} We prove the two claims of the theorem in the 
two respective parts below.
 
\textbf{1.} Firstly, 
$r_{2} \textbf{y} = \textbf{y} r_{2} = 0$, 
$r_{2} \textbf{x} = \textbf{x} r_{2}$, $\left( -r_{2}\right) ^{2}=-r_{2}$.
Therefore,
$r = - r_2$ 
is a central idempotent, and 
$
\mathcal{R} = r \mathcal{R} \oplus (1-r) \mathcal{R}
 = r \mathbb{Z}\left\langle \textbf{x}\right\rangle \oplus (1-r) \mathcal{R}$
where $(1-r)\mathcal{R} \cong M_n\left( \mathbb{Z}\right)\,$,
and $ r \mathbb{Z}\left\langle \textbf{x}\right\rangle \cong \mathbb{Z} C_n$.

\textbf{2.} We construct a finite set, call it
$\mathcal{S}$,
such that every element of $S(H)$ may be written as an integral linear combination
of the elements of $\mathcal{S}$. 

Multiply the relation $r_{2,n}(x,y) = 0$ by $y$ on the left:
\begin{equation} \label{m8}
y^2 + \sum_{i = 1}^{n-1} yx^{-i}yx^{i}-y=0 \Longrightarrow
y^{2}-y=-\sum_{h\in H}yx^{h}yx^{-h}\text{.}
\end{equation}
Therefore, for $k\in H$, we have 
\begin{equation}\label{m9}
\left( y-y^2 \right) x^{k}y = \sum_{h\in H}yx^{h}yx^{-h}x^{k}y =  yx^{k}y^{2}.
\end{equation}
Multiply the relation $r_{2,n}(x,y)  = 0$ by $y$ on the right:
\begin{equation}\label{m7} 
y^{2}+\sum_{i=1}^{n-1}x^{-i} y x^{i} y - y=0 \Longrightarrow
y^{2}-y=-\sum_{h\in H} x^{-h} y x^{h} y.
\end{equation}\label{m15}
It follows that
\begin{equation}\label{m15}
y^{3}-y^{2}=-y\sum_{h\in H}x^{-h}yx^{h}y=0.
\end{equation}
Therefore, multiplying (\ref{m9}) by $y$ on the left yields
\begin{equation}\label{m10}
y^2 x^{k} y^{2} = - \left( y^{3}-y^2\right) x^{k}y = 0 \text{.}
\end{equation}
Let $k \in H$, then equating the right-hand sides of (\ref{m8}) and (\ref{m7}) 
gives us
\begin{equation}\label{m12}
yx^{k}yx^{-k} = -\sum_{k \neq h\in H}yx^{h}yx^{-h}+\sum_{h\in
H}x^{-h}yx^{h}y
\end{equation}
Next, multiplying (\ref{m12}) by $y x^l$ on the left and
by $x^k$ on the right yields
\begin{equation}\label{m14}
yx^{l} y x^{k} y = -\sum_{h\in H,h\not=k}yx^{l} yx^{h}yx^{-h+k} +
\sum_{h\in H} y x^{l} x^{-h} y x^{h} y x^{k} = y^{2}x^{l}yx^{k}.
\end{equation}
We conclude that every word in $x$ and $y$ may be rewritten in such a way that
the following conditions are satisfied:
\begin{enumerate}
\item $x$ occurs finitely many times with exponent between $0, \ldots, n-1$,
because one of the relation in (\ref{Grigdream}) is $x^n = 1$.
\item Powers of $y$ may occur as subwords at most most twice because of
(\ref{m14}). 
\item $y$ occurs with exponent between $0,1,2$
because $y^3 = y^2$ by (\ref{m15}).

\end{enumerate}
Stated another way, every element in $S(H)$ may be written
as $\mathbb{Z}$-linear combination of
the words of the form 
$x^{\alpha_1} y^{\beta_1} x^{\alpha_2} y^{\beta_2} x^{\alpha_3}$,
where $\alpha_1, \alpha_2, \alpha_3 \in \{0, \ldots, n-1 \}$ and
$\beta_1, \beta_2 \in \{0,1,2 \}$.
\end{proof}

\subsection{Magnus-type ring extension of $M_n(\mathbb{Z})$}
In the proof of Theorem \ref{Said3} below, we 
introduce an analog of the Magnus Embedding
from Magnus \cite{Magnus} (see Lemma on p.~764 of \cite{Magnus}).
 
\begin{theorem}\label{Said3}
The ring $\mathbb{Z} \{ x,y\}$ has a quotient $\mathcal{R} = \mathcal{R}_n$ such that
\begin{enumerate}
\item $\mathcal{R}$ is an over-ring of $M_n(\mathbb{Z})$.
\item Under the natural epimorphism $\mathcal{R} \twoheadrightarrow M_n(\mathbb{Z})$, 
the images of the ideals generated by $r_{1n}, s_1, \ldots, s_n$ form a direct sum.
\end{enumerate}  
\end{theorem}

\begin{proof} The proof consists of finding a ring $\mathcal{R}$ 
such that

\smallskip

1.  $\mathcal{R}$ is generated by two elements $\mathbf{x}, \mathbf{y}$ 
together with $1_{\mathcal{R}}$.

\smallskip

2. Let $\mathcal{R}_1  = \mathcal{R} r_{1,n}( \mathbf{x}) \mathcal{R}$, 
$\mathcal{S}_i = \mathcal{R} s_i( \mathbf{x}, \mathbf{y}) \mathcal{R}$ for
$1 \leq i \leq n-1$, and $\mathcal{S}_0 = \mathcal{R} s_0( \mathbf{y}) \mathcal{R}$.
Then $\mathcal{R}_1 \cap \mathcal{S}_0 = \{ 0_{\mathcal{R}} \}$ and
$\mathcal{S}_0 = \mathcal{S}_1 \oplus \ldots \oplus \mathcal{S}_{n-1}$. 

\medskip

Put $M = M_n(\mathbb{Z})$ and 
consider the ring 
$\mathcal{M} = \left( 
\begin{array}{cc}
M & 0 \\ 
\xi  M\oplus \eta 
 M & \mathbb{Z} 
\end{array}
\right) $ 
where 
$\xi$ and $\eta $ 
are independent variables commuting with each other and with every matrix from
$M$. 
{\em Let $\mathcal{R}$ be the subring of $\mathcal{M}$
generated by the matrices} 
\[
\textbf{I} = \left( 
\begin{array}{cc}
I & 0 \\ 
0 & 1
\end{array}
\right) ,
\,
\textbf{X}=\left( 
\begin{array}{cc}
X & 0 \\ 
\xi & 1
\end{array}
\right) ,
\,
\textbf{Y}=\left( 
\begin{array}{cc}
Y & 0 \\ 
\eta & 0
\end{array}
\right) .
\]
Then the projection on the top left corner is a ring epimorphism 
$ \mathcal{R} \twoheadrightarrow M$, by Theorem \ref{Grig}.  Define the polynomials
$q_{0}(t)=0$ and 
$q_{i}(t)=1+t+...+t^{i-1}$, $i\geq 1$. Then
\[
\mathbf{X}^{i} =\left( 
\begin{array}{cc}
X^{i} & 0 \\ 
\xi q_{i}(X) & 1
\end{array}
\right),~
\mathbf{X}^{-i}=\left( 
\begin{array}{cc}
X^{-i} & 0 \\ 
-\xi q_{i}(X) X^{-i} & 1
\end{array}
\right),~
\mathbf{Y} \mathbf{X}^{i} \mathbf{Y} =\left( 
\begin{array}{cc}
0 & 0 \\ 
\eta X^{i} Y  & 0
\end{array}
\right),
\]
\[
\mathbf{X}^{i} \mathbf{Y} =\left( 
\begin{array}{cc}
X^{i} Y & 0 \\ 
\xi q_{i}(X) Y+\eta & 0
\end{array}
\right),~
\mathbf{X}^{-i} \mathbf{Y}=\left( 
\begin{array}{cc}
X^{-i} Y & 0 \\ 
-\xi q_{i}(X) X^{-i}Y + \eta & 0
\end{array}
\right),
\]
\[
\mathbf{X}^{-i} \mathbf{Y} \mathbf{X}^{i} =\left( 
\begin{array}{cc}
X^{-i} Y X^{i} & 0 \\ 
-\xi q_{i}(X) X^{-i} Y X^{i} + \eta X^{i}  & 0
\end{array}
\right) =
\left( 
\begin{array}{cc}
X^{-i} Y X^{i} & 0 \\ 
\xi q_{i}(X) X^{-i} Y X^{i}  +\eta X^{i}) & 
0
\end{array}
\right) \text{.}
\]

For the remainder of the proof, let $r_1 = r_{1,n}(\mathbf{X})$, 
$s_j = s_j(\mathbf{X},\mathbf{Y})$, and $1\leq i\leq n-1$. Then
\[
r_{1}  =\left( 
\begin{array}{cc}
0 & 0 \\ 
\xi q_{n}(X) & 0
\end{array}
\right) ,
s_{0} =\left( 
\begin{array}{cc}
0 & 0 \\ 
\eta \left( Y -1\right) & 0
\end{array}
\right) \text{,} 
~
s_{i} =\left( 
\begin{array}{cc}
0 & 0 \\ 
\eta  X^{i}Y & 0
\end{array}
\right) \text{.}
\]
Therefore,
\[
\mathcal{R}_{1}=\left( 
\begin{array}{cc}
0 & 0 \\ 
\xi q_{n}(X)M & 0
\end{array}
\right) ,
~
\mathcal{S}_{0}=\left( 
\begin{array}{cc}
0 & 0 \\ 
\eta \left( Y-1\right) M & 0
\end{array}
\right),~
\mathcal{S}_{i}=\left( 
\begin{array}{cc}
0 & 0 \\ 
\eta X^{i}Y \,M & 0
\end{array}
\right) \text{.}
\]
We see that $\mathcal{R}_{1}\cap \mathcal{S}_0= \{ 0  \}$.
The significance of this fact will become apparent from
the following claim that will finally prove the theorem.

\textbf{Claim.} {\em The sum $ \sum _{i = 1}^{n-1}\mathcal{S}_{i}$ is direct and 
equals $\mathcal{S}_{0}$.}

We argue as follows. An element $u_{0}$ in $\mathcal{S}_{0}$ has the form 
$u_{0}=\left( 
\begin{array}{cc}
0 & 0 \\ 
\eta T_{0} & 0
\end{array}
\right)$
where
$T_{0}=\left( Y-1\right) M_{0}$ 
for some
$M_{0}=\left( 
\begin{array}{c}
M_{01} \\ 
\vdots \\ 
M_{0n}
\end{array}
\right) \in M
$.
Therefore, $T_{0}=-\left( 
\begin{array}{c}
0 \\ 
M_{02} \\ 
\vdots \\ 
M_{0n}
\end{array}
\right)$.
An element $u_{i}$ in $\mathcal{S}_{i}$ $\left( 1\leq i\leq n-1\right) $ has
the form
$
u_{i}=\left( 
\begin{array}{cc}
0 & 0 \\ 
\eta T_{i} & 0
\end{array}
\right)
$
where
$T_{i}=X^{i}Y M_{i}$
for some
$M_{i}=\left( 
\begin{array}{c}
M_{i1} \\ 
\vdots \\ 
M_{in}
\end{array}
\right) \in M\text{.}
$
Then the $i+1$st row of $T_i$ is $M_{i1}\,$,
the other rows being  zero. Therefore, 
$\sum_{i = 1}^{n-1} T_{i} = \left( 
\begin{array}{c}
0 \\ 
M_{11} \\ 
\vdots \\ 
M_{n-1,1}
\end{array}
\right)$ is of the same form as $T_0 \,$, and hence
$\sum_{i = 1}^{n-1} T_{i} \in \mathcal{S}_0$. We conclude that $\sum_{i = 1}^{n-1}T_{i}=0$ if and only if $
M_{i1}=0$ for all $i \in \{1, \ldots, n-1\}$.
\end{proof}

\subsection{$M_n(\mathbb{Z})$ as a quotient of rings without identity}\label{Fantasy}

To motivate this discussion, let $\mathcal{R} = \mathbb{Z}  \{ e_{11}, \ldots, e_{nn} \} $ be 
a free non-associative ring without identity.
Let $\mathcal{I}$ be the ideal of $\mathcal{R}$ generated by the elements
$e_{ij} e_{kl} - \delta_{jk} e_{il}$. Then the quotient ring
$\mathcal{R} / \mathcal{I}$ is isomorphic to $M_n(\mathbb{Z})$.

Another way to present $M_n(\mathbb{Z})$ as a quotient of a ring without identity
is to modify Presentation (\ref{Grigdream}) to obtain $M_n(\mathbb{Z})$   
as a quotient of the integral semigroup ring 
$\mathbb{Z}[FS( x,y )]$.
This yields the following
\begin{theorem}\label{NoIdentity}
Let $X = \sum_{i=1}^n E_{i,i+1}$ and $Y = E_{11}$. Then the map
\[
f: \mathbb{Z}[FS ( x,y )] \to M_n(\mathbb{Z}),~ x \mapsto X , ~ y \mapsto Y
\]
is a ring epimorphism with kernel generated by the $n+2$ elements
\begin{equation}\label{GenSet}
x^{n+1}-x,~ yx^{n}-y, ~ 
-x^{n} + \sum_{i = 0}^{n-1}x^{n-i}yx^{i},~
yx^{j}y, ~1\leq j\leq n-1\text{.}
\end{equation} 
\end{theorem}

\begin{proof}
Put $\mathcal{R} = \mathbb{Z}[FS ( x,y )]$. 
{\em All computations in this paragraph will be done modulo $\mathcal{I} = \text{Ker}(f)$.}
We firstly observe that $x^{n-1}\left( x^{n+1}-x\right) = 0$ yields $x^{2n}=x^{n}$. Therefore
$y^2  = y\left( \sum_{i = 0}^{n-1}x^{i}yx^{n-i}\right) = yx^{n} = y$, so that
$y = y^2 = \left( \sum_{i = 0}^{n-1}x^{i}yx^{n-i}\right) y = x^{n}y $.
Therefore, $z=x^{n}$ an identity element.

It remains to show that ideal $\mathcal{I}_0$ generated by the elements (\ref{GenSet}) 
equals $\mathcal{I}$. 
Firstly, $\mathcal{I}_0 \subseteq \mathcal{I}$ because the corresponding relations 
are satisfied by $X$ and $Y$. On the other hand, the computations in the previous
paragraph show that the ring $\mathcal{R} / \mathcal{I}_0$ 
is generated by the $n^2$ elements $x^i + \mathcal{I}_0$, 
$x^i y x^j + \mathcal{I}_0$, $1 \leq i,j \leq n$. Since 
$\dim _{\mathbb{Z}} M_n(\mathbb{Z}) = n^2$, it follows that $\mathcal{I} = \mathcal{I}_0$.
\end{proof}

\subsection{Linear representations of matrix rings}
We prove below that $4$
relations in $X$ and $Y$ are sufficient to describe $M_{n}\left( \mathbb{Z}\right)$
in the context of matrix rings.

\begin{theorem}\label{Said1}
Let $\mathcal{D}$ be a commutative domain of characteristic either zero or 
at least $m+1$,
over which every finitely generated 
projective module is free.
Let
$\mathcal{S} $ 
be a subring of $M_{m}\left( \mathcal{D}\right)$ generated 
by some nonzero $X_1$ and $ Y_1 $
such that
\begin{equation}\label{1eq}
X_1^{n+1} = X_1, ~ Y_1 X_1^n = Y_1, ~ Y_1^2 = Y_1,~
\sum_{i = 0}^{n-1} X_1^{n-i} Y_1 X_1^{i} = X_1^{n}
\text{.}
\end{equation}
Then the trace $k$ of $Y_1$ is a positive integer, and
there exist $B \in GL_m(\mathcal{D})$ such that, putting
$r = m - kn$, we have
\begin{equation*}
B^{-1}X_{1}B=\left( 
\begin{array}{cc}
I_{k}\otimes X & 0_{k\times r} \\ 
0_{r\times k} & 0_{r\times r}
\end{array}
\right)~\, \text{and}~ \,
B^{-1}Y_{1}B=\left( 
\begin{array}{cc}
I_{k}\otimes Y & 0_{k\times r} \\ 
0_{r\times k} & 0_{r\times r}
\end{array}
\right) \text{.}
\end{equation*}
\end{theorem}

An exposition of commutative domains over which every finitely generated
projective module is free can be found in Lam \cite{Lam2}.

\begin{proof}[Proof of Theorem \ref{Said1}]

In this paragraph, the ring $\mathcal{S}$ will be embedded
into a ring smaller than $M_m(\mathcal{D})$.
Since $X_{1}^{n}$ is an idempotent, we decompose $\mathcal{D}^{m}$ as
the direct sum of the image $\mathcal{P}$ and the kernel $\mathcal{N}$, i.e.~ 
$\mathcal{D}^{m}= \mathcal{P} \oplus \mathcal{Z} $ where
\begin{enumerate} 
\item $\mathcal{P}$ and $\mathcal{Z}$ have 
$\mathcal{D}$-ranks  $q$ and $r$, respectively.
\item $X_{1}^{n}|_{P}=I_{q}$ and $X_{1}^{n}|_{Z}=0_{r}$. 
\end{enumerate}
We observe from (\ref{1eq}) that $\mathcal{P}$ and $\mathcal{Z}$ are 
$\mathcal{S}$-invariant and $\mathcal{S}|_{\mathcal{Z}}=0_{r}$. Choose some free
generating sets for $\mathcal{P}$ and $\mathcal{Z}$. 
Then with respect to these sets, $X_{1}$  and $Y_1$ are
represented by  the matrices $\left( 
\begin{array}{cc}
X_{2} & 0 \\ 
0 & 0
\end{array}
\right) $ and  
$\left( 
\begin{array}{cc}
Y_{2} & 0 \\ 
0 & 0
\end{array}
\right) $, respectively.
Furthermore, the matrices $X_2$ and $Y_2$ satisfy the following relations
\begin{equation}\label{11eq}
r_{1,n}(X_2, Y_2) = r_{2,n}(X_2, Y_2) = s_0( Y_2) = 0.
\end{equation}

\smallskip

Let $k = \text{tr}(Y_2)$. Then (\ref{11eq}) yield 
\begin{equation}\label{2eq}
q = \text{tr}(I_{q}) = \text{tr} \left( 
\sum_{i = 0}^{n-1} X_{2}^{i} Y_{2} X_{2}^{n-i} \right) =
\sum_{i = 0}^{n-1} \text{tr} \left( Y_{2} X_{2}^{n-i} X_{2}^{i}  \right) = nk.
\end{equation}
$\mathcal{P}$ decomposes 
with respect to 
the idempotent $Y_{2}$ as a direct sum of
the image $\mathcal{U}$ and the kernel $\mathcal{V}$.
The restriction maps $Y_{2}|_{\mathcal{U}}$
and $Y_{2}|_{\mathcal{V}}$
are the identity and zero maps, respectively. Therefore
\begin{equation}\label{3eq}
k = \text{tr}(Y_2) = \text{tr} \left( Y_{2}|_{\mathcal{U}} \right) 
+ \text{tr} \left( Y_{2}|_{\mathcal{V}} \right) = \text{tr} \left( Y_{2}|_{\mathcal{U}} \right)
 = \text{tr} \left( id_{\mathcal{U}} \right) \text{.}
\end{equation}
In particular, $k$ is an integer. 

Let
\[ 
\widehat{\mathcal{U}} = \sum_{i=0}^{n-1} X_{2}^i\left( \mathcal{U}\right) \text{.}
\] 
Then (\ref{1eq}) implies that $\widehat{\mathcal{U}}$
is an  
$\mathcal{S}$-module. In addition, $Y_{2}|_{\mathcal{V}} = 0 $
yields $Y_{2}|_{\mathcal{P} / \widehat{\mathcal{U}}} = 0 $. In turn,
(\ref{1eq}) implies $X_{2}|_{\mathcal{P} / \widehat{\mathcal{U}}} = 0 $,
which amounts to the identity map acting as zero on 
$\mathcal{P} / \widehat{\mathcal{U}}$.  
Therefore 
$  \mathcal{P} =\widehat{\mathcal{U}}$. The sum 
$\sum_{i=0}^{n-1} X_{2}^i\left( \mathcal{U}\right)$
is direct because by passing to the field of fractions $\mathcal{F}$ of
$\mathcal{D}$, we have 
$\mathcal{F}^q = \sum_{i=0}^{n-1} X_{2}^i\left(\mathcal{F} \otimes_{\mathcal{D}} \mathcal{U}\right)$.
By (\ref{3eq}), this sum is a sum of $n$ linear spaces of dimension $k$, and we know
from (\ref{2eq}) that $\dim_{\mathcal{F}} \mathcal{F}^q = nk$.
Therefore 
\[
\mathcal{D}^{q}= \bigoplus_{i=0}^{n-1} X_{2}^i\left( \mathcal{U}\right).
\]

Let $\mathcal{B} = \{s_1, \ldots, s_k \}$ be a free $\mathcal{D}$-basis
of $\mathcal{U}$. Then 
$\widehat{\mathcal{B}} = \bigcup_{i=0}^{n-1} X_2^{i}(\mathcal{B})$
is a free $\mathcal{D}$-basis of $\mathcal{P}$.
Hence, $X_2$ may be represented with respect to 
$\widehat{\mathcal{B}}$ 
by an $n$-by-$n$ block matrix $\left( X_{ij}\right) $
with $k$-by-$k$ blocks,
where $X_{ij}=0$ unless $i=j+1$, and $X_{j+1,j} = I_k$ for $1 \leq j \leq n-1$.
Similarly, 
$Y_{2}=\left( Y_{ij}\right) $ where $Y_{11}=I_{k}$ and $Y_{ij} = 0$
for $i \neq 1$ because $Y_2 |_{U}$ is the identity map,
and $Y_2 |_{\mathcal{P} / \mathcal{U}}$ is the zero map.
Since 
$I_q = X_{2}^{n}=X_{1,n} \otimes I_{n}$, we arrive at $
X_{1,n}=I_{k}$. Therefore, $X_2$ is represented in the basis
$\widehat{\mathcal{B}}$
by the permutation matrix $X \otimes I_k$
in block form. It remains to observe that
from $r_{2,n}\left( X_2, Y_2 \right) = I_q$ it follows that
$Y_{1j} = 0$ for $2 \leq j \leq n$. Consequently
$Y_2$ is represented with respect to $\widehat{\mathcal{B}}$
by the matrix $Y \otimes I_k$.
\end{proof}

\begin{corollary}\label{Said2}
Let $\mathcal{D}$ be a commutative domain of characteristic either zero
or at least $n+1$. Then 
the automorphism group of the ring $M_n \left( \mathcal{D}\right) $ is generated 
by the automorphism group $Aut(\mathcal{D})$ of the ring $\mathcal{D}$, and 
by the projective general linear group $ PGL_n \left(\mathcal{D} \right)$,
where
\begin{enumerate}
\item $Aut(\mathcal{D})$ acts on $M_n \left( \mathcal{D}\right) $
by acting on each entry of a matrix.
\item $ PGL_n \left(\mathcal{D} \right)$ acts on $M_n \left( \mathcal{D}\right) $
by conjugation.
\end{enumerate}
\end{corollary}

\begin{proof}
Any automorphism $\sigma$ of the ring $M_n \left( \mathcal{D}\right) $
leaves the center invariant. In other words, there exist 
$\alpha \in Aut(\mathcal{D})$ such that for every
$a \in \mathcal{D}$, we have 
$\sigma\left( a \sum_{i=1}^n E_{ii} \right) = \alpha(a) \sum_{i=1}^n E_{ii}$.

Next we consider $ \beta = \alpha^{-1} \sigma $, which is a 
$\mathcal{D}$-algebra automorphism of $M_n(\mathcal{D})$. 
Then the pair 
$\left( \beta X, \beta Y \right)$ 
satisfies the relations of (\ref{Dnepr}). Therefore,  by Theorem \ref{Said1}
there exists $U\in M_{n}\left( \mathcal{D}\right) $ which conjugates 
$\beta X$ to $X$ and $\beta  Y$ to $Y$. The conjugations by
$U$ and $-U$ produce identical results, and there are no further such identifications.
Therefore the automorphism group of the $\mathcal{D}$-algebra $ M_n \left( \mathcal{D}\right)$
is isomorphic to $PGL_n \left(\mathcal{D} \right)$.  
\end{proof}

The result of Corollary \ref{Said2} is not new. More general results
are contained Rosenberg and Zelinsky \cite{Rosenberg}. In particular, that
paper shows that Corollary \ref{Said2} is false, for example, for Dedekind
domains with class number at least $2$.

\smallskip

We will need the following theorem of G.~Higman \cite{Higman}.

\begin{theorem}[G.~Higman's Theorem]\label{HigmanTh} 
The unit group $\mathcal{U}$ of the integral group ring of a finite 
Abelian group $\mathcal{G}$ is given by 
$\mathcal{U} = \pm \, \mathcal{G} \times \mathcal{F}$, where
$\mathcal{F}$ is a free Abelian group of rank
\begin{equation}\label{HigmanEq}
\frac{1}{2} \left(\# \mathcal{G} + t_2  - 2l + 1 \right).
\end{equation}
Here $t_2$ is the number of elements of $\mathcal{G}$ of order $2$,
and $l$ is the number of cyclic subgroups of $\mathcal{G}$. 
\end{theorem}
\noindent{By analyzing some elementary inequalities, it follows that $\mathcal{F} = \{ 0 \}$ if and only
if $n = 2,3,4,6$.}

\begin{theorem}\label{Said-Bogdan} 
The set 
$\mathcal{Y} = \{Y_1 \in M_n(\mathbb{Z}) \mid  Y_1^2 = Y_1, ~ r_{2,n}(X,Y_1) = 0\}$ 
has the property that the pair $(X, Y_1)$  
satisfies all relations of (\ref{Dnepr}), and
all $Y_1$ have trace $1$.
If $n = 2,3,4,6$ then 
$Y_{1} = E_{ii}$ for some $i$.
Otherwise, $\mathcal{Y}$
is
infinite, and if $Y_1 \neq E_{ii}$
then it has both positive and negative entries.

Any $Y_1$ is of the form $\left(c_i d_j \right)$ for some
integers $c_i, d_j$ such that
the matrices $\text{circ} \left(c_1, \ldots, c_n  \right)$ and 
$\text{circ} \left(d_1, \ldots, d_n  \right)$ are mutually inverse.
Any $Y_1$ is conjugate to $Y$ by an integral circulant matrix
with determinant $\pm 1$.
\end{theorem}

\begin{proof}
Let $Y_1 = (y_{ij})$. Then $r_{2,n} \left(X, Y_1  \right) = 0$
implies
\begin{equation}\label{StJoe}
\sum_{k=0}^n y_{i+k, \, j+k} = \delta_{ij}. 
\end{equation}
These formulas
prove the claim about the possible signs of entries of $Y_1$.

Applying the trace to $r_{2,n} (X,Y_1) = 0$ implies

\begin{equation}\label{2eq}
n = \text{tr}(I_{n}) = \text{tr} \left( 
\sum_{i = 0}^{n-1} X^{i} Y_{1} X^{n-i} \right) =
\sum_{i = 0}^{n-1} \text{tr} \left( Y_{1} X^{n-i} X^{i}  \right) = n \, tr(Y_1).
\end{equation}

$\mathbb{Z}^n$ decomposes with respect to the idempotent $Y_1$ as a direct sum
of the image $\mathcal{I}$ and kernel $\mathcal{K}$. Therefore
\begin{equation}\label{3eq}
1 = \text{tr}(Y_1) = \text{tr} \left( Y_{1}|_{\mathcal{I}} \right) 
+ \text{tr} \left( Y_{1}|_{\mathcal{K}} \right) = \text{tr} \left( Y_{1}|_{\mathcal{I}} \right)
 = \text{tr} \left( id_{\mathcal{I}} \right) \text{.}
\end{equation}
Therefore, $Y_1$ is a rank $1$ projection.
The image
of $Y_1$ is
an Abelian group is generated by some $\left(d_1, \ldots, d_n  \right) \in \mathbb{Z}^n$.
It follows that on the standard basis $e_1 = (1,0, \ldots, 0), \ldots, e_n = (0, \ldots, 0,1)$ the action of $Y_1$
is described by $Y_1 e_i = c_i d_1 + \ldots + c_i d_n$ for some integer $c_i$. Therefore
$Y_1 = \left( c_i d_j \right)$. Next, from $r_{2,n} \left( X, Y_1 \right) = 0$
we conclude that $\sum_{k = 0}^{n-1} c_{i+k} d_{j+k} = \delta_{ij}\,$, which
is the same as saying that the matrices $\text{circ}  \left( c_1, \ldots, c_n \right)$
and $\text{circ}  \left( d_1, \ldots, d_n \right)$ are mutually inverse.

Now, going back to (\ref{StJoe}), we see that the relations $Y_1 X^k Y_1 = 0$ follow
from the relations $r_{1,n}(X) = r_{2,n}(X, Y_1) = 0$. Indeed,
$\left( X^k  Y_1 \right)_{ij} = c_{i+k}d_j$. Therefore 
$
\left(Y_1 X^k  Y_1 \right)_{ij} = c_i \left( \sum_{u=1}^n d_u c_{u+k} \right) d_j =
c_i \delta_{kn} d_j = 0$.
It follows that $(X,Y_1) \in G_n(\mathbb{Z})$ by Theorem \ref{Grig} and because
all proper quotients of the ring $M_n(\mathbb{Z})$ are finite. 

In the cases of $n = 2,3,4,6$ the group $\mathcal{U}(\mathbb{Z} \langle X  \rangle)$
consists precisely of $2n$ matrices $\pm E_{ii}$.        
\end{proof}

Theorem \ref{Said-Bogdan} may be strengthened as follows. If
all entries of $X_1 \in M_n(\mathbb{Z})$ are nonnegative,
and $X_1^n = I_n\,$, then in each row of $X_1$ there exactly one positive entry, and it equals
$1$. 
We will prove this assertion in $2$
steps.

\smallskip

\textbf{1.} Suppose that in each row of $X_1$ there is exactly one
nonzero entry. Then from $\det X_1 = \pm 1$ it follows that $X_1$
is of the required form.

\smallskip

\textbf{2.} 
Suppose that $X_1 = (x_{ij})$ has a row with at least $2$ positive entries
$x_{ij}$ and $x_{ij'}$. The $i$th column of $X_1$ contains
a nonzero entry $x_{mi}$. We conclude that the matrix $X_1^2 = (t_{kl})$
has the property that $t_{mj}, \, t_{mj'} > 0$. Similarly,
any positive power of $X_1$ has at least two positive entries in some
row. We obtain a contradiction, however, by considering $X_1^n = I_n$.

\medskip

We remark that G.~Higman's Theorem \ref{HigmanTh},
when applied to a cyclic group of order $n$,
may be restated in terms of solutions of the following Diophantine equations:
\begin{equation}\label{Circ}
\det \text{circ}(x_1, \ldots, x_n) = \pm 1.
\end{equation}
Unfortunately, there appears to be no efficient algorithm to find solutions of (\ref{Circ}).
Computer experiments with (\ref{Circ}) eventually
led us to Theorem \ref{Said-Bogdan}.

\subsection{Presentations of direct sums of matrix rings over $\mathbb{Q}$ and $\mathbb{Z}$}

Our next result shows that the ring $M_n(\mathbb{Z})$ has infinitely many presentations.
We obtain, as a consequence, the presentations for several types of direct sums of matrix
rings. We do not write down these presentations explicitly based on the following
reason.  
If $\mathcal{I}$
and $\mathcal{J}$ are ideals of a ring $\mathcal{R}$ such that
$\mathcal{I} + \mathcal{J} = \mathcal{R}$, then 
$\mathcal{I} \cap \mathcal{J} = \mathcal{I} \mathcal{J} + \mathcal{J} \mathcal{I}$.
Therefore, if the ideals $\mathcal{I}$ and $\mathcal{J}$ are
generated by explicitly given $i$ and $j$ elements, respectively, then
$\mathcal{I} \cap \mathcal{J}$ is generated by
at most $ 2ij$ explicitly given elements. 
\begin{theorem}\label{Modular}
The ring $\mathbb{Z} \{ x,y \}$ has an infinite family of ideals
$\{ \mathcal{I}_n(m) \}_{m \in \mathbb{Z}}$ defined by
\[
\mathcal{I}_n (m) = 
\left(r_{1,n}(x, mx + y),~ r_{2,n}(x, mx + y),~ s_j,~ 1 \leq j \leq n-1  \right),~
\mathcal{I}_n = \mathcal{I}_n (0).
\]
This family of ideals has the following properties:
\begin{enumerate} 
\item\label{chu1} $\mathbb{Z} \{ x,y \} / \mathcal{I}_n(m) 
\cong M_n(\mathbb{Z})$ for any integer $m$.

\item\label{chu2} If $S$ is a finite subset of $\mathbb{Z}$ and
$\mathcal{R} = \mathbb{Z} \{x,y \} / \bigcap_{s \in S} \mathcal{I}_n(s)$, then \\
$\mathcal{R} \otimes_{\mathbb{Z}} \mathbb{Q} \cong M_n(\mathbb{Q})^S$.

\item\label{chu2.5} If $n_1, \ldots, n_k$ are
pairwise relatively prime integers, then \\
$
\mathbb{Z} \{ x,y \} / \bigcap_{j=1}^k \mathcal{I}_{n_j}
\cong \bigoplus_{j=1}^{k} M_{n_j} (\mathbb{Z})$.

\item\label{chu3} If $|k-l| \geq 2$, then even though
$\mathbb{Z} \{x,y \} / \mathcal{I}_n(k) \cap \mathcal{I}_n(l) 
\ncong M_n(\mathbb{Z})^2 $,
it embeds as a subring of finite index.

\item\label{chu4} Define the map $^t: \mathbb{Z} \{x,y \} \to \mathbb{Z} \{x,y \}$
by $f(x,y)^t = f(y,x)$, then
$
\mathbb{Z} \{x,y \} / \mathcal{I}_n ^t \cap \mathcal{I}_n \cap \mathcal{I}_n(1)
\cong M_n(\mathbb{Z})^3  $.

\item\label{chu5}
$
\mathbb{Z} \{x,y \} / \mathcal{I}_2 \cap \mathcal{I}_2 ^t 
\cap \mathcal{I}_2(1) \cap \mathcal{I}_2(1)^t
\cong M_2(\mathbb{Z})^4  $. 
\end{enumerate}
\end{theorem}

\begin{proof}
We find it convenient to introduce a family of ring automorphisms 
$\{\varphi_m \}_{m \in \mathbb{Z}}$
of $\mathbb{Z} \{ x,y \}$ given by $\varphi_m (x) = x$ and
$\varphi_m (y) = mx + y$. 
Then 
$\mathcal{I}_n(m) = \varphi_m \left(\mathcal{I}_n \right)$.
Theorem \ref{Grig} tells is that $\mathbb{Z} \{ x,y \} / \mathcal{I}_n(m) \cong M_n(\mathbb{Z})$.

We will show that all ideals $\mathcal{I}_n(m)$ are different.
Suppose that this is false, so that 
$
\mathcal{I}_n(k) = \mathcal{I}_n(l)
$
for some $k \neq l$. 
Then
\[
\mathcal{I}_n = \varphi_{-k} \left( \mathcal{I}_n(k) \right) = \mathcal{I}_n(l-k)=
\mathcal{I}_n(s),~\text{where}~ 0 \neq a = l-k.
\]
Therefore,
\begin{itemize}
\item $r_{2,n}(x,ax+y) = a \sum_{i=0}^{n-1} x^{n+1} + r_{2,n}(x,y) \equiv
nsx \, \left(\text{mod}\, \mathcal{I}_n \right)$  

\item $r_{2,n}(x,ax+y) \in \mathcal{I}_n$

\item $x$ is invertible modulo $\mathcal{I}_n$
\end{itemize}
imply that $na \in \mathcal{I}_n$.
Therefore,
$\{ 0 \} = na \left(\mathbb{Z} \{x,y \} / \mathcal{I}_n \right)
\cong M_n(\mathbb{Z})$, a contradiction.
The argument above, together with Chinese Remainder Theorem, proves Parts \ref{chu1} 
and \ref{chu2}.

\smallskip

To prove Part \ref{chu2.5}, we will show that 
$\mathcal{I}_{ij} = \mathcal{I}_{n_i} + \mathcal{I}_{n_j} 
= \mathbb{Z}\{ x,y \}$ when $i \neq j$, and the computations will be done
modulo $\mathcal{I}_{ij}$. From
$x^{n_i} = x^{n_j} = 1$ it follows by Euclid's Algorithm that
$x = x^{\text{gcd}(n_i, n_j)} = 1$
$\Longrightarrow$ $0 = yxy = y^2 = y$ $\Longrightarrow$
$0 = r_{2,n_i}(x,y) = r_{2,n_i}(x,0) = -1$.

\smallskip

Next we prove 
Part \ref{chu3}. We observe that the restriction of the maps $\varphi_m$ 
to $\mathbb{Z}$ is the identity map. Therefore,
\[
\left(\mathcal{I}_n(k) + \mathcal{I}_n(l) \right) \cap \mathbb{Z} = 
\left(\mathcal{I}_n + \mathcal{I}_n(k - l) \right) \cap \mathbb{Z} 
\equiv \mathcal{I}_n \cap \mathbb{Z} 
\equiv \{0\}(\text{mod}~ k-l ),
\] 
so that 
$\mathcal{I}_n(k) \cap \mathcal{I}_n(k) \cap \mathbb{Z} \neq \mathbb{Z}$.
To sum up, 
$
n(k-l) \mathbb{Z} \subseteq 
\left(\mathcal{I}_n(k) + \mathcal{I}_n(l) \right) \cap \mathbb{Z}
\subsetneqq \mathbb{Z}
$
yields 
$\mathcal{I}_n(k) + \mathcal{I}_n(l) \neq \mathbb{Z} \{ x,y \}$.

\smallskip

We prove Part \ref{chu4} by showing that the sum of any two of the three ideals
$\mathcal{I}_n ^t, \, \mathcal{I}_n, \, \mathcal{I}_n(1) $ is $\mathbb{Z} \{x,y \}$.

\textbf{1.} We claim that 
$\mathcal{J} = \mathcal{I}_n + \mathcal{I}_n(1)  = \mathbb{Z} \{x,y \}$. All
computations here are done modulo $\mathcal{J}$. 
We observe that
\begin{eqnarray*}
0 = s_j (x,x+y) = (x+y) x^j (x+y) &=& \\
x^{j+2} + x^{j+1} y + y x^{j+1} + s_j (x,y) & = &
x^{j+2} + x^{j+1} y + y x^{j+1} = v_j(x,y). 
\end{eqnarray*}
Therefore
$0 = v_{n-2}(x,y) = x^n + x^{n-1} y + y x^{n-1} = 1+ x^{n-1} y + y x^{n-1}$,
which we multiply by $y$ on the left and by $x$ on the right:
\begin{eqnarray}\label{by1}
0 = yx + (y x^{n-1} y)x + y^2 x^{n} = yx + y.
\end{eqnarray}
Likewise,
\begin{multline}\label{by2}
0 = x v_{n-1}(x,y) y = x(1+ x^{n-1} y + y x^{n-1})y =
xy + y^2  + s_{n-1}(x,y)  = \\
 xy + y. 
\end{multline}
One consequence of (\ref{by1}) and (\ref{by2}) is 
$ 0 = - s_1(x,y) = - yxy = y^2 = y$, 
and we conclude that
$ 0 = r_{2,n}(x,y) = -1$.

\textbf{2.} We claim that $\mathcal{K} = \mathcal{I}_n ^t + \mathcal{I}_n = 
\mathbb{Z} \{ x,y \}$. 
All computations here are done modulo $\mathcal{K}$.
From $x^2 = x$ and $x^n = 1$ we conclude that
$0 = yxy = y^2 = y$, and therefore $0 = 0^n = y^n = 1$.

\textbf{3.} We claim that 
$\mathcal{L} = \mathcal{I}_n ^t + \mathcal{I}_n(1) = 
\mathbb{Z} \{ x,y \}$. All computations here are done modulo $\mathcal{L}$.
As above, $x = 1$. Then $0 = xyx = y$, and therefore, as above,
$0 = 0^n = y^n = 1$.

\textbf{4.} We claim that $\mathcal{M} = \mathcal{I}_n + \mathcal{I}_n(1) = 
\mathbb{Z} \{ x,y \}$. All computations here are done modulo $\mathcal{M}$:
$0 = y(x+y)y = yxy + y^2 = y$, then $0 = r_{2,n}(x,y) = -1 $.

\smallskip

It remains to prove Part \ref{chu5} of the theorem. In view
of the arguments of Part \ref{chu4}, it remains to show that
$\mathcal{N} = \mathcal{I}_2 (1) + \mathcal{I}_2 (1)^t = \mathbb{Z} \{ x,y \}$,
and as usual, all the necessary computations will be done
modulo $\mathcal{N}$.
\begin{equation}\label{Zpush1}
 x + y = (x+y)^2 = x^2 + xy + yx + y^2 = 2 + xy + yx. 
\end{equation}
\begin{equation}\label{Zpush2}
0 = (x+y) x (x+y) y = (x^3 + x^2 y + y x^2 + yxy)y = 
xy + 2 + yx. 
\end{equation}
The right-hand sides of (\ref{Zpush1}) and (\ref{Zpush2})
are equal, hence $x+y = 0$.
Therefore, $0 = r_{2,2}(x,x+y) = r_{2,2}(x,0) = -1$. 
\end{proof} 

\smallskip

While by Part \ref{chu3} of Theorem \ref{Modular}, it is already impossible
to obtain $M_2(\mathbb{Z})^5$ as a quotient of $\mathbb{Z} \{ x,y \}$
by intersecting the ideals $\mathcal{I}_2(m)$ and $\mathcal{I}_2(m)^t$,
we wonder whether there exists an infinite family 
$\{ \mathcal{T}_m \}_{m \geq 1}$ of ideals in 
$\mathbb{Z}\{ x,y \}$ such that 
$\mathbb{Z} \{x,y \} / \bigcap_{m=1}^k \mathcal{T}_m \cong M_2(\mathbb{Z})^k $.
One possible obstacle to overcome here would be to use
various subgroups of 
$\text{Aut} \left( \mathbb{Z} \{ x,y \}  \right)$ to create new ideals from
$(x^2-1, \, y + xyx - 1, \, yxy)$, and then determine their
interdependence.
In the proof of Theorem \ref{Modular}, for example,
we have used an infinite cyclic subgroup and a subgroup
of order $2$. There is a similar question about an arbitrary 
$M_n(\mathbb{Z})$ as well.


\begin{thebibliography}{9}
 


\bibitem{Burnside1} W.~Burnside,
{\em On the condition of reducibility of a group of linear substitutions.}
Proc.~London.~Math.~Soc.~3 (1905) 430-434.


\bibitem{Cohen} H.~Cohen,
{\em A course in computational algebraic number theory.} 
Graduate Texts in Mathematics, 138.
Springer-Verlag, Berlin, 1993.

\bibitem{Curtis} C.W.~Curtis and I.~Reiner, {\em Representation theory of 
finite groups and associative algebras.} Pure and Applied Mathematics, 
Vol. XI Interscience Publishers, a division of John Wiley \& Sons, New York-London 1962.



\bibitem{Faith} C.~Faith,
{\em Algebra I Rings, Modules and Categories},
Springer-Verlag, 1981.

\bibitem{Frohlich}  A.~Fr\"{o}hlich and M.J.~Taylor, {\em Algebraic number theory.} 
Cambridge Studies in Advanced Mathematics, 27. Cambridge University Press, Cambridge, 1993.


\bibitem{Hensley} D.~Hensley, private communication.


\bibitem{Higman} G.~Higman,
{\em The units of group-rings.}  
Proc. London Math. Soc. (2)  46, (1940). 231--248.



\bibitem{Lam2}  T.Y.~Lam,
{\em Serre's conjecture.}
Lecture Notes in Mathematics, Vol. 635.
Springer-Verlag, Berlin-New York, 1978.



\bibitem{Lam} T.Y.~Lam, {\em A first course in noncommutative rings.} 
Second edition. Graduate Texts in Mathematics, 131. Springer-Verlag, New York, 2001.
  
 


\bibitem{Long} W.E.~Longstaff, {\em Burnside's theorem: irreducible pairs of transformations.}  
Linear Algebra Appl.  382  (2004), 247--269.

\bibitem{Magnus} W.~Magnus, 
{\em On a theorem of Marshall Hall.}  
Ann. of Math. (2)  40,  (1939). 764--768.


 



\bibitem{Rosenberg} A.~Rosenberg and D.~Zelinsky, 
{\em Automorphisms of separable algebras.}
Pacific J. Math. 11 1961 1109--1117.



\bibitem{Saltman} D.J.~Saltman, private communication.

\bibitem{Sehgal}  S.K.~Sehgal,
{\em Topics in group rings.}
Monographs and Textbooks in Pure and Applied Math., 50.
Marcel Dekker, Inc., New York, 1978.


\end{thebibliography}
\end{document}